\documentclass[10pt,a4paper,american, reqno]{amsart}
\usepackage{xcolor}
\usepackage{babel}
\usepackage{amsthm}
\usepackage{amsbsy}
\usepackage{amstext}
\usepackage{amssymb}
\usepackage{esint}
\usepackage[colorlinks,pdfpagelabels,pdfstartview = FitH,bookmarksopen = true,bookmarksnumbered = true,linkcolor = blue,plainpages = false,hypertexnames = false,citecolor = red, pagebackref=true]{hyperref}
\usepackage{breakurl}
\usepackage{enumitem}
\usepackage{graphicx}

\makeatletter

\numberwithin{equation}{section}
\numberwithin{figure}{section}

\usepackage{babel}
\usepackage{amsfonts}
\usepackage{amsthm}
\usepackage{mathcomp}
\usepackage{dsfont}
\usepackage{latexsym}

\def\Xint#1{\mathchoice
    {\XXint\displaystyle\textstyle{#1}}%
    {\XXint\textstyle\scriptstyle{#1}}%
    {\XXint\scriptstyle\scriptscriptstyle{#1}}%
    {\XXint\scriptscriptstyle\scriptscriptstyle{#1}}%
    \!\int}

\newtheorem{theo}{Theorem}[section]

\newtheorem{prop}[theo]{Proposition}
\newtheorem{lem}[theo]{Lemma}
\newtheorem{definition}[theo]{Definition}
\newtheorem{remark}[theo]{Remark}

\theoremstyle{definition}

\newtheorem*{lem*}{Lemma}

\newtheorem*{cor*}{Corollary}
\newtheorem*{theo*}{Theorem}

\DeclareMathOperator{\Div}{div}
\DeclareMathOperator{\diam}{diam}

\newcommand{\N}{\ensuremath{\mathbb{N}}}
\newcommand{\R}{\ensuremath{\mathbb{R}}}

\newcommand{\BIGOP}[1]{\mathop{\mathchoice%
{\raise-0.22em\hbox{\huge $#1$}}%
{\raise-0.05em\hbox{\Large $#1$}}{\hbox{\large $#1$}}{#1}}}

\newcommand{\mint}{- \mskip-19,5mu \int}

\def\mwint_#1{\mathchoice
 {\mathop{\vrule width 6pt height 3 pt depth -2.5pt
        \kern -8.5pt \intop}\nolimits_{\kern -3pt #1}}%
 {\mathop{\vrule width 6pt height 3 pt depth -2.6pt
                  \kern -6pt \intop}\nolimits_{#1}}%
 {\mathop{\vrule width 6pt height 3 pt depth -2.6pt
                  \kern -6pt \intop}\nolimits_{#1}}%
 {\mathop{\vrule width 6pt height 3 pt depth -2.6pt
                  \kern -6pt \intop}\nolimits_{#1}}}

\newcommand{\dx}{\mathrm{d}x}

\newcommand{\dt}{\mathrm{d}t}
\newcommand{\ds}{\mathrm{d}s}
\newcommand{\dtau}{\mathrm{d}\tau}
\newcommand{\dsigma}{\mathrm{d}\sigma}

\renewcommand{\epsilon}{\varepsilon}

\def\LL{L}
\def\Lip{{\rm Lip}^x}

\numberwithin{equation}{section}

\def\Xint#1{\mathchoice
    {\XXint\displaystyle\textstyle{#1}}%
    {\XXint\textstyle\scriptstyle{#1}}%
    {\XXint\scriptstyle\scriptscriptstyle{#1}}%
    {\XXint\scriptscriptstyle\scriptscriptstyle{#1}}%
    \!\int}
\def\XXint#1#2#3{\setbox0=\hbox{$#1{#2#3}{\int}$}
    \vcenter{\hbox{$#2#3$}}\kern-0.5\wd0}

\def\dashint{\Xint{\raise4pt\hbox to7pt{\hrulefill}}}

\def\Xiint#1{\mathchoice
    {\XXiint\displaystyle\textstyle{#1}}%
    {\XXiint\textstyle\scriptstyle{#1}}%
    {\XXiint\scriptstyle\scriptscriptstyle{#1}}%
    {\XXiint\scriptscriptstyle\scriptscriptstyle{#1}}%
    \!\iint}
\def\XXiint#1#2#3{\setbox0=\hbox{$#1{#2#3}{\iint}$}
    \vcenter{\hbox{$#2#3$}}\kern-0.5\wd0}
\def\biint{\Xiint{-\!-}}

\allowdisplaybreaks

\makeatother

\author[V. B\"ogelein]{Verena B\"{o}gelein}
\address{Verena B\"ogelein\\
Fachbereich Mathematik, Universit\"at Salzburg\\
Hellbrunner Str. 34, 5020 Salzburg, Austria}
\email{verena.boegelein@plus.ac.at}

\author[F. Duzaar]{Frank Duzaar}
\address{Frank Duzaar\\
Fachbereich Mathematik, Universit\"at Salzburg\\
Hellbrunner Str. 34, 5020 Salzburg, Austria}
\email{frankjohannes.duzaar@plus.ac.at}

\author[G. Treu]{Giulia Treu}
\address{Giulia Treu\\
Dipratimento di Matematica 'Tullio Levi-Civita', Universita' di Padova\\
Via Trieste 63, 35121 Padova, Italy }
\email{giulia.treu@unipd.it}

\begin{document}
\title[Parabolic PDEs with Dynamic Data under a Bounded Slope Condition]{Parabolic PDEs with Dynamic Data under a Bounded Slope Condition}

\begin{abstract} 
We establish the existence of Lipschitz-continuous solutions to the Cauchy–Dirichlet problem for a class of evolutionary partial differential equations of the form
\begin{equation*}
    \partial_tu-\Div_x \nabla_\xi f(\nabla u)=0
\end{equation*}
in a space-time cylinder $\Omega_T=\Omega\times (0,T)$, subject to time-dependent boundary data
$g\colon \partial_{\mathcal{P}}\Omega_T\to \R$
prescribed on the parabolic boundary. The main novelty in our analysis is a time-dependent version of the classical bounded slope condition, 
imposed on the boundary data $g$
along the lateral boundary
$\partial\Omega\times (0,T)$.  More precisely, we require that for each fixed $t\in [0,T)$, the graph of $g(\cdot ,t)$ over $\partial\Omega$
admits supporting hyperplanes with slopes that may vary in time but remain uniformly bounded. The key to handling time-dependent data lies in constructing more flexible upper and lower barriers.
\end{abstract}

\subjclass[2010]{35A01, 35K61, 35K86, 49J40}
\keywords{Existence, parabolic equations, bounded slope condition, Lipschitz solutions}

\maketitle
\tableofcontents

\section{Introduction}

Throughout this paper, let $\Omega \subset \mathbb{R}^n$, $n \geq 2$, denote a bounded, open, and convex set, and let $f \colon \mathbb{R}^n \to \mathbb{R}$ be a convex integrand. A classical theorem of Haar \cite{Haar:1927} ensures that, for prescribed boundary values $u_o \colon \partial \Omega \to \mathbb{R}$ satisfying the \emph{bounded slope condition} (BSC), there exists a Lipschitz continuous minimizer $u \colon \Omega \to \mathbb{R}$ of the variational functional
\begin{equation}
\mathcal{F}(v) := \int_{\Omega} f(Dv) \,\dx, \label{eq:F-energy} 
\end{equation}
subject to the boundary condition $u = u_o$ on $\partial \Omega$; see also \cite{Hartmann-Nirenberg, Hartmann-Stampacchia, Miranda:1965, Stampacchia}. The construction has become standard and appears in modern textbooks on the calculus of variations; see, e.g., \cite[Chapter~1]{Giusti.2003}.
More recent developments concerning the existence of Lipschitz minimizers under the BSC may be found in \cite{Bousquet:2007, Bousquet:2010, Cellina:2001, Clarke:2005, Mariconda-Treu:2002, Mariconda-Treu:2007, Mariconda-Treu:2009, Mascolo-Schianchi:1983}. The standard argument hinges on a comparison principle in which the affine functions furnished by the BSC serve as barriers, together with the translation invariance of minimizers and a priori gradient bounds. The successful application of this strategy depends sensitively on the structure of the integrand, and in particular breaks down when lower-order terms are present. In such cases, affine functions no longer provide suitable barriers, and the method must be adapted by employing a more flexible class of barrier functions; see \cite{Fiaschi-Treu:2012, Giannetti-Treu:2022, Giannetti-Treu:2025}.

Surprisingly, a time-dependent counterpart to this semi-classical theory has remained an open problem for some time. Instead, sophisticated parabolic methods -- such as Galerkin approximations, monotone operators, and nonlinear semigroup theory -- have led to various existence results. However, the construction of Lipschitz continuous solutions to evolutionary equations associated with general convex integrands $f$
has remained elusive without additional assumptions on the growth of the integrand. The paper \cite{BDMS_bd-slope} marks the first Haar-Rado-type result for the Cauchy-Dirichlet problem
\begin{equation*}
\left\{
\begin{array}{cc}
	\partial_tu-\Div_x \nabla_\xi f(\nabla u)=0, & \quad\mbox{in $\Omega_T=\Omega\times (0,T)$,}\\[4pt]
	u=u_o, & \quad \mbox{on $\partial_{\mathcal{P}}\Omega_T$,}
\end{array}
\right.
\end{equation*}
which guarantees the existence of a unique classical solution $u\in C^0([0,T];L^2(\Omega)) \cap L^\infty (\Omega_T)$ with Lipschitz continuous spatial gradient
$\nabla u\in L^\infty(\Omega_T,\R^n)$.
The only conditions required are the convexity of the integrand $f\colon\R^n\to \R$, and the BSC imposed on the initial and lateral boundary datum $u_o\in W^{1,\infty}(\Omega)$. 
A related result for linear growth functionals can be found in \cite{Hardt-Zhou}; see also \cite{Zhou}. The theory has subsequently been further generalized to data $u_o$
satisfying a one-sided BSC, as shown in \cite{Boeg-Sta, Rai-Sil-Sta, Sta}. A common feature of all these parabolic papers is the use of the BSC to construct affine barriers.

At this juncture, the natural question arises as to whether a semi-classical theory can be developed for time-dependent boundary values. Specifically, one might ask whether it is possible to construct Lipschitz solutions to the Cauchy-Dirichlet problem \eqref{eq:strong-Cauchy-Dirichlet-prob} with a boundary datum 
$g\colon \partial_{\mathcal{P}}\Omega_T\to\R$
depending on time. As we will discuss later, the proof strategy of \cite{BDMS_bd-slope} cannot, in principle, be applied to time-dependent boundary conditions. The barrier construction employed in \cite{BDMS_bd-slope} necessitates that the boundary values $g$
be independent of time.

\subsection{Lipschitz solutions}
In the parabolic setting, the formulation of Lipschitz continuous variational solutions involves certain function spaces that can be interpreted as the parabolic analogue of Lipschitz continuous functions from the stationary setting. Specifically, we use the identification between the space of Lipschitz continuous functions, $\Lip(\Omega) = C^{0,1}(\Omega)$, and the Sobolev space $W^{1,\infty}(\Omega)$ to define the function space consisting of those bounded functions $v \in L^\infty(\Omega_T)$ having a bounded spatial gradient $\nabla v \in L^\infty(\Omega_T, \mathbb{R}^n)$, i.e.,
\[
	\Lip(\Omega_T)
	:=
	\big\{ v\in L^{\infty}(\Omega_T)\colon
	\nabla v\in L^\infty(\Omega_T,\R^n) \big\}.
\]
For a given $\LL \in (0, \infty)$, we define the subclass $\Lip(\Omega_T, \LL)$ as
\[
	\Lip(\Omega_T, \LL)
	:=
	\big\{ v\in \Lip (\Omega_T)\colon
	\|\nabla v\|_{L^\infty(\Omega_T,\R^n)}\le\LL \big\}.
\]
As a consequence of the identification
of  $\Lip(\Omega_T)$ with $L^\infty (0,T, W^{1,\infty}(\Omega))$, for almost every $t \in (0,T)$, the restriction $v(t) := v(\cdot, t)$ of $v \in \Lip(\Omega_T)$ to the time slice $\Omega \times \{t\}$ belongs to $W^{1,\infty}(\Omega)$. 
This allows the classical definition of the trace on appropriate time slices. The space of Lipschitz continuous functions $\Lip_0(\Omega_T)$ with zero lateral boundary values is also well-defined, enabling the formulation of the Cauchy-Dirichlet problem in the parabolic function space $g + \Lip_0(\Omega_T)$ for some $g\colon\Omega_T\to\R$. This space consists of those functions $v \in \Lip(\Omega_T)$ such that for almost every time slice $t \in (0,T)$ the restriction $v(t)$ of $v$ to $\Omega \times \{t\}$ belongs to $g(t) + W^{1,\infty}_0(\Omega) \equiv g(t) + C^{0,1}_0(\Omega)$. 
For the boundary values $g$, we assume 
\begin{align}\label{ass:g} 
\text{ $g \in \Lip(\Omega_T)$ with $\partial_t g \in L^2(\Omega_T)$, and $g(0) := g_o \in W^{1,\infty}(\Omega)$.} 
\end{align}
Instead of $g+\Lip_0(\Omega_T)$ we write $\Lip_g(\Omega_T)$ for those $v\in \Lip (\Omega_T)$ coinciding with
$g$ on the lateral boundary $\partial\Omega\times (0,T)$.

For the definition of the notion of {\em variational solution} to the parabolic Cauchy-Dirichlet problem
\begin{equation}\label{eq:strong-Cauchy-Dirichlet-prob}
\left\{
\begin{array}{cc}
	\partial_tu-\Div_x \nabla_\xi f(\nabla u)=0, & \quad\mbox{in $\Omega_T=\Omega\times (0,T)$,}\\[4pt]
	u=g, & \quad \mbox{on $\partial_{\mathcal{P}}\Omega_T$}
\end{array}
\right.
\end{equation}
we follow the approach of Lichnewsky and Temam \cite{Lichnewsky.1978}, originally introduced in the context of the time-dependent parametric minimal surface equation. This idea leads naturally to the notion of variational solutions used in the following. 

\begin{definition}[Variational Solution]\label{def:var-sol} \upshape Assume that the Cauchy-Dirichlet datum
$g$ fulfills hypothesis \eqref{ass:g}.  A map $u\in \Lip_g(\Omega_T)$ is called \emph{variational solution on $\Omega_{T}$} to the Cauchy-Dirichlet
problem \eqref{eq:strong-Cauchy-Dirichlet-prob} if and only if the
variational inequality 
\begin{align}
	\iint_{\Omega_\tau}f(\nabla u)\,\dx\dt
	& \le
	\iint_{\Omega_\tau}\big[\partial_{t}v(v-u)
	+
	f(\nabla v)\big]\,\dx\dt\nonumber \\
 	& \quad+
	\tfrac{1}{2}\|v(0)-g_{o}\|_{L^{2}(\Omega)}^{2}
	-
	\tfrac{1}{2}\|(v-u)(\tau)\|_{L^{2}(\Omega)}^{2}\label{eq:variational-inequality}
\end{align}
holds true, for  a.e.~$\tau\in [0,T]$ and for any $v\in \Lip_{g}(\Omega_T)$
with $\partial_{t}v\in L^{2}(\Omega_{T})$ and $v(0)\in L^2(\Omega)$. 
\end{definition}
Note that all terms in \eqref{eq:variational-inequality} are well-defined, since $u \in \Lip_g(\Omega_T)$ implies $u \in L^\infty(0,T; L^2(\Omega))$. By assumption \eqref{ass:g}, the function $g$ is admissible as a test function in the variational inequality \eqref{eq:variational-inequality}. This permits testing with $v = g$, which yields certain energy estimates. As a consequence, one can conclude that variational solutions satisfy the initial condition $u(\cdot, 0)=g_o$ in the usual $L^2(\Omega)$ sense, as demonstrated in Lemma~\ref{lem:initial}. The concept of a variational solution enables the use of techniques from the calculus of variations. Under certain conditions, it can be demonstrated that a variational solution also satisfies the properties of a weak solution; see Theorem~\ref{thm:reg-gen} below.

As in the classical theory of variational problems, we begin by constructing variational solutions subject to a \emph{gradient constraint}. In this regard, we introduce a notion of variational solution that incorporates such a constraint. Assuming that for some $\LL > 0$, we have
\begin{equation}\label{Def:L}
	\LL>\max\big\{ \|\nabla g_o\|_{L^\infty(\Omega,\R^n)}, \|\nabla g\|_{L^\infty(\Omega_T,\R^n)}\big\},
\end{equation}
the subclass $\Lip_g(\Omega_T, \LL)$ is non-empty, and the following definition of variational solutions with gradient constraints is well-posed.
\begin{definition}[Variational solutions, gradient constraint]\label{def:var-sol-rest}
\upshape 
Let $g$ be as in \eqref{ass:g} and $\LL >0$ as in \eqref{Def:L}.
A map $u\in \Lip_g(\Omega_T,\LL)$ 
is called  \emph{variational solution of the gradient constrained Cauchy-Dirichlet problem} \eqref{eq:strong-Cauchy-Dirichlet-prob}
in $\Lip_g(\Omega_T,\LL)$ if and only if the
variational inequality 
\begin{align}
	\iint_{\Omega_\tau}f(\nabla u)\dx\dt
	& \le
	\iint_{\Omega_\tau}\big[\partial_{t}v(v-u)
	+
	f(\nabla v)\big]\dx\dt\nonumber \\
 	& \quad+
	\tfrac{1}{2}\|v(0)-g_{o}\|_{L^{2}(\Omega)}^{2}
	-
	\tfrac{1}{2}\|(v-u)(\tau)\|_{L^{2}(\Omega)}^{2}\label{eq:variational-inequality-rest}
\end{align}
holds true, for a.e.~$\tau\in [0,T]$ and  for any  $v\in \Lip_g(\Omega_T,\LL)$ 
with $\partial_{t}v\in L^{2}(\Omega_{T})$ and $v(0)\in L^2(\Omega)$.  \hfill$\Box$
\end{definition}

\subsection{The main results}

Our main result in this paper is the existence of Lipschitz continuous solutions (with respect to the spatial variable) to the Cauchy-Dirichlet problem \eqref{eq:strong-Cauchy-Dirichlet-prob} with time-dependent boundary values, under the assumption of the bounded slope condition. Before stating the main result, we first outline the \emph{assumptions on the data}.

\begin{enumerate}[label=(A\arabic{*}), ref=(A\arabic{*}), topsep=7pt]
    \item\label{A1}
    Let $\Omega\subset\R^n$ be $R$-uniformly convex for some $R>0$ in the sense of Definition~\ref{def:unif-conv}; 
    \item\label{A2} 
    Let $f\colon \R^n\rightarrow\R$ be convex, and, outside the unit ball 
    $B_1=B_1(0)$, of class $C^{1,1}$ and  uniformly convex.
    Specifically,  there exists $\epsilon\in(0,1]$, such that
    \begin{equation}\label{hess}
        D^2f(\xi) (\zeta, \zeta)\ge \epsilon|\zeta|^2
        \qquad\mbox{for a.e. $\xi\in\R^n\setminus B_1$ and any $\zeta\in\R^n$;}
  \end{equation}
\item\label{A3} 
Let $g\in W^{1,\infty}(\Omega_T)$ be such that its restriction to the lateral boundary, $g|_{\partial\Omega\times [0,T]}$, 
satisfies the time-dependent bounded slope condition $t-{\rm BSC}_Q$ in the sense of Definition \ref{def:bdslope}, for some constant $Q>0$.
In addition, assume that for each $x_o\in \partial\Omega$, the associated functions $w_{x_o}^\pm$	from the bounded slope condition are Lipschitz continuous in time, that is,
$$
w_{x_o}^\pm\in W^{1,\infty}([0,T],\R^n),
$$
and that  
$$
\mathsf Q:=
\sup_{x_o\in\partial\Omega}\|(w_{x_o}^\pm)'\|_{L^\infty([0,T],\R^n)}<\infty.
$$
\end{enumerate}

Then the following existence result holds.

\begin{theo}[Existence of Lipschitz solutions]\label{thm:existence-var-sol}
Suppose that assumptions \textnormal{\ref{A1}}--\textnormal{\ref{A3}} are satisfied. Then there exists a unique variational solution 
to the Cauchy--Dirichlet problem \textnormal{\eqref{eq:strong-Cauchy-Dirichlet-prob}} in the sense of Definition~\ref{def:var-sol}\, which 
satisfies the gradient bound
$$
	\|\nabla u\|_{L^\infty (\Omega_T,\R^n)}
	\le 
	C,
$$
where $C$ depends on $n, \epsilon, R, \diam(\Omega), f, \nabla f, \|D^2f\|_{L^\infty(\R^n\setminus B_1)}, Q, [g]_{0,1;\Omega_T}$,  and  $\mathsf Q$.
\end{theo}

We emphasize that the variational solution is unique even when the integrand $f$ is convex but not strictly convex; see \cite[Lemma~3.3]{BDS}. This includes the case of integrands with flat regions. Moreover, if $f\in C^1$, then the variational solution constructed above enjoys additional regularity.

\begin{theo}[Regularity of solutions] 
\label{thm:reg-gen}
Suppose that assumptions \textnormal{\ref{A1}}--\textnormal{\ref{A3}} are satisfied, and that $f\in C^1$. Then the variational solution $u$
obtained in Theorem~\ref{thm:existence-var-sol} is a weak solution to the Cauchy--Dirichlet problem~\eqref{eq:strong-Cauchy-Dirichlet-prob}, and satisfies $u\in C^{0;1,\frac12}(\Omega_T)$.
\end{theo}

\subsection{Novelty and key technical tools}\label{sec:novelty}
The proof proceeds in two steps. The first establishes the existence of variational solutions for time-dependent boundary data, under a uniform Lipschitz bound on the admissible class. This reduces to a variational inequality with a gradient constraint. The construction follows De Giorgi’s minimizing movement scheme, cf.~\cite{Ambrosio-Gigli-Savare}, which is particularly suited to incorporating time-dependent Dirichlet data on the lateral boundary into the class of admissible competitors, cf.~\cite{Bogelein.2014b, BDS, Schatzler}. An alternative derivation, based on the method of weighted energy dissipation, cf.~\cite{Akagi-Stefanelli, Bogelein.2013, Bogelein.2014, Mielke-Stefanelli:2011}, draws on De Giorgi’s variational framework for nonlinear evolution equations, including applications to certain nonlinear wave equations; see \cite{DeGiorgi.1996, Ilmanen:1994, Serra.2012, Stefanelli}.
This may be viewed as the first step in the spirit of Haar’s \cite{Haar:1927} approach to constructing Lipschitz minimizers of the parametric area functional—a perspective that anticipated modern developments in geometric analysis and regularity theory for minimal surfaces.

The solutions obtained are variational solutions in the sense of Lichnewsky and Temam \cite{Lichnewsky.1978}; that is, they solve a variational inequality subject to a gradient constraint of the form $|\nabla u|\le L$, where 
$L>0$ is a fixed constant. As a result, the class of admissible variations is restricted to those that preserve this structural bound. In particular, the presence of the constraint limits perturbations to directions that remain within the admissible set. It is therefore necessary to establish that the variational solution lies strictly inside the admissible class; only then can arbitrary variations be performed, and the solution identified as the sole admissible one satisfying the variational inequality in the unconstrained class.

In the classical setting, this difficulty was addressed by Hilbert \cite{Hilbert} through the introduction of the bounded slope condition, later refined by Haar \cite{Haar:1927} in the construction of Lipschitz minimizers of the parametric area functional. By restricting attention to graphs with uniformly bounded gradient, compactness is restored and direct methods become applicable, obviating the need for unconstrained variations. The bounded slope condition furnishes barriers required for comparison and maximum principles, and renders the explicit gradient constraint superfluous by enforcing uniform control through the prescribed bound.

The use of affine barriers derived from the bounded slope condition fundamentally relies on temporal constancy of the Dirichlet boundary data. This constraint becomes evident upon inspecting the argument in \cite[Theorem~1.2]{BDMS_bd-slope}, where for a fixed point $x_o\in\partial \Omega$
the function $g(x_o,t)+w_{x_o}^-(t)\cdot(x-x_o)$ is employed as a lower barrier and must act as a sub-solution to \eqref{eq:strong-Cauchy-Dirichlet-prob}$_1$. Differentiation in time yields the constraint
$$
    \partial_t\big(g(x_o,t)+w_{x_o}^-(t)\cdot(x - x_o)\big)
    =
    \partial_t g(x_o,t)+ \partial_t w_{x_o}^-(t)\cdot(x-x_o)
    \le 
    0
$$
for any $(x,t)\in\Omega_T$. Taking the limit $x\to x_o$ within $\Omega$ implies $\partial_t g(x_o,t)\le 0$ for any $t\in(0,T)$.  A symmetric argument involving an  upper barrier of the form $g(x_o,t)+w_{x_o}^+(t)\cdot(x-x_o)$  leads  to the inequality $\partial_t g(x_o,t)\ge 0$. One is thus forced to conclude that $g$ is stationary on
the lateral boundary $\partial\Omega\times (0,T)$.  In the presence of genuinely time-dependent boundary values, such affine constructions therefore fail to apply. While the bounded slope condition retains its structural role, the analysis in the presence of  time-dependent boundary data necessitates more flexible constructions: the use of affine sub- or super-solutions, though formally admissible, inherently restricts the boundary values to be time-independent and is therefore incompatible with the temporal variability intrinsic to the problem.


Instead of employing affine barriers, we implement a construction based on the convex conjugate $f^\ast$, following an observation due to Cellina~\cite{Cellina:2007}. The central idea exploits the fact that, when the integrand $f$ is of class $C^2$, a function of the form
\begin{equation*}
    v(x,t)
    =
    \frac{n}{\alpha}f^*\Big(\frac{\alpha}{n}\big(x-y(t)\big)\Big)-c(t),
\end{equation*}
with $\alpha\in\R\setminus\{0\}$, $y\colon[0,T]\to\R^n$, and $c\colon[0,T]\to\R$ satisfies 
\begin{equation*}
    \Div_x\nabla_\xi f\big(\nabla v(x,t)\big)=\alpha.
\end{equation*}
In particular, if $\alpha >0$, then a choice of $\alpha$, $y$, and $c$ such that $\alpha$ dominates the time derivative 
$\partial_tv$ at every point 
$(x,t)\in \Omega_T$
would ensure that $v$
is a sub-solution to the parabolic equation~\eqref{eq:strong-Cauchy-Dirichlet-prob}, and thus a potential candidate for a barrier from below.
However, being a subsolution is only one of the requirements for a function to qualify as a barrier from below. In addition, it must be compatible with the boundary data in the sense that, for some fixed boundary point $x_o\in\partial\Omega$, one has
$ v(x_o,t)=g(x_o,t)$ for any $t\in [0,T)$, and $v\le g$
on $\partial_{\mathcal{P}}\Omega_T$.
At this stage, the principal difficulty lies in ensuring that the real parameter 
$\alpha$, together with the functions $y$ and $c$, can indeed be chosen so that 
$v$ constitutes a lower barrier. This is realized by considering, for each $t\in [0,T)$, the set 
\begin{equation*}
    \widetilde\Omega_t
    :=
    \bigg\{x\in\R^n: \underbrace{\frac{n}{\alpha}f^*\Big(\frac{\alpha}{n}\big(x-y(t)\big)\Big)-c(t)}_{=\, v(x,t)}\le g(x_o,t)+\widetilde w_{x_o}^{-}(t)\cdot (x-x_o)\bigg\},
\end{equation*}
where $\widetilde w_{x_o}^{-}$ is a suitable modification of the affine function $w_{x_o}^{-}$, chosen to satisfy 
$g(x_o,t)+\widetilde w_{x_o}^{-}(t)\cdot (x-x_o)\le g(x,t)$ for all $(x,t)\in\Omega_T$. 
We then carefully select
$y(t)$ and $c(t)$ such,
on the one hand,
\begin{equation*}
    x_o\in\partial \widetilde\Omega_t
    \quad\mbox{for all $t\in [0,T)$,}
\end{equation*}
and, on the other hand,
\begin{equation*}
\Omega\subset\widetilde\Omega_t
    \quad\mbox{for all $t\in [0,T)$.}
\end{equation*}
By the definition of the sets $\widetilde\Omega_t$, this implies $v(x_o,t)=g(x_o,t)$ for all $t\in [0,T)$
and $v\le g$ in $\Omega_T$. Next, we establish an upper bound for $|\partial_t v|$ in $\Omega_T$ that is independent of $\alpha$. This enables us to choose $\alpha$
sufficiently large to guarantee that $v$
is a subsolution. Finally, we compute an upper bound on the spatial gradient $\nabla v$, thereby demonstrating that the constructed lower barrier is Lipschitz continuous. Realizing this approach would already be delicate for $C^2$-integrands. The weaker assumption \ref{A2} on $f$
requires even greater care to balance the parameters involved and to avoid circular dependencies.

\subsection{Explicit construction of a lower barrier in a specific example}
To illustrate the construction of the lower barriers summarized in the preceding subsection and described in detail in \S\,\ref{sec:remove}, we provide a simple example in which these barriers are explicitly constructed. Thereby, we use the notation introduced in \S\,\ref{sec:novelty} and \S\,\ref{sec:remove}. Let us consider the Cauchy-Dirichlet problem \eqref{eq:strong-Cauchy-Dirichlet-prob}, where $\Omega =B_1(0)\subset\R^2$, 
$f(\xi)=\frac12 |\xi|^2$ with convex conjugate $f^\ast (\eta)=\frac12 |\eta|^2$,
and boundary data  
$g(x,t)=x_1\cos t+x_2\sin t =(\cos t ,\sin t)\cdot (x_1,x_2)$.
We observe that, in this case, the boundary values are affine in $x=(x_1,x_2)$
for each fixed $t$, and they satisfy the 
$t-{\rm BSC}_Q$ condition
with constant $Q=1$. In particular, we note that the $t-{\rm BSC}_Q$  is satisfied for any $x_o\in \partial B_1(0)$ with the choice  $w_{x_o}^-(t)=(\cos t,\sin t)$. In the setting of this example, the function defined in \eqref{def:barrier} takes the form
\begin{equation*}
    v(x,t)=\tfrac{\alpha}{4}|x-y(t)|^2-c(t),
\end{equation*}
where $\alpha >0$ is a parameter, and 
$y(t)\in\R^2$, $c(t)\in \R$
are functions of time. We fix the point 
$x_o=(-1,0)\in \partial B_1(0)$ and aim to show that one can choose $\alpha$, $y$, and $c$ such that the following conditions are satisfied:
\begin{enumerate}[label=(C\arabic{*}), ref=(C\arabic{*}), topsep=7pt]
\item\label{C1} $v(x_o,t)=g(x_o,t)$ for all $t\in[0,T)$;
\item\label{C2} $v(x,t)\le g(x,t)$ for all $(x,t)\in\partial_{\mathcal{P}}\big(B_1(0)\times[0,T)\big)$;
\item\label{C3} $v$ is a sub-solution of
\begin{equation}\label{eq:example}
\partial_t v-\Delta v =0 \qquad\text{\rm in } B_1(0)\times[0,T).
\end{equation}
\end{enumerate}
For $\alpha>0$ we now choose $y(t)$ and $c(t)$ such that \ref{C1} and \ref{C2} are satisfied. To this end, we define, cf.~\S\,\ref{sec:remove},
\begin{align*}
    \widetilde{\Omega}_t
    &=
    \big\{x\in\R^2: \tfrac{\alpha}{4}|x-y(t)|^2-c(t)\le g(x_o,t)+w_{x_o}^-(t)\cdot(x-x_o)\big\} \\
    &=
    \big\{x\in\R^2: \tfrac{\alpha}{4}|x-y(t)|^2-c(t)\le x_1 \cos t+x_2\sin t\big\}.
\end{align*}
Note that $g(x_o,t)+w_{x_o}^-(t)\cdot(x-x_o)= g(x,t)$ for any $(x,t)\in B_1(0)\times[0,T)$ and hence we can take $\widetilde w_{x_o}^-=w_{x_o}^-$ in 
Lemma~\ref{lem:tBSC-G}. 
Next, we choose
$$
y(t)=(y_1(t),y_2(t))=\tfrac{2}{\alpha}\left(1-\cos t,-\sin t\right)
$$ 
and 
$$
c(t)= \tfrac{\alpha}{4}+1+\tfrac{2}{\alpha}\left(1-\cos t\right),
$$ 
and demonstrate that
$x_o=(-1,0)\in\partial\widetilde{\Omega}_t$ and $B(0,1)\subset\widetilde{\Omega}_t$.
Through direct calculations, we obtain
\[
    \widetilde{\Omega}_t
    = B_{1+\frac{2}{\alpha}}\big(\tfrac{2}{\alpha},0\big)\quad \mbox{ for every $t\in[0,T)$,}  
\] 
and hence $\widetilde{\Omega}_t$ satisfies the desired properties.
To verify \ref{C3}, i.e.~that $v(x,t)$
is a subsolution of \eqref{eq:example} for appropriate values of $\alpha >0$,
we proceed as follows.
By the properties of $f$ and $f^*$ stated in Proposition \ref{propertyf}, it immediately follows that $ \Delta v=\Div_x\nabla_\xi f(\nabla v)=\alpha$ in $B_1(0)\times[0,T)$. By a straightforward computation, we obtain for the time derivative of $v$ that
\begin{align*}
\partial_t v(x,t)
&=-x_1\sin t+x_2\cos t .  
\end{align*}
We can conclude that $v$ is a subsolution of \eqref{eq:example} in the set $B_1(0)\times[0,T)$ for any $\alpha\ge 1$. Moreover, we estimate that $|\nabla v|\le 2+\frac\alpha2$ in $B_1(0)\times[0,T)$.

In the next Figure \ref{fig:1} we represent the above construction for the fixed value $\alpha=2$, at time $t=0$ on the left and at time $t=\frac{\pi}{4}$ on the right. To be more precise, in the picture on the left there are the affine function involved in the Bounded Slope Condition  for the point $x_o=(-1,0)$ and the corresponding barrier; the continuous line is the boundary condition, while the dashed line is the trace of the barrier function on  $\partial\widetilde\Omega_0$. The picture on the right has to be interpreted analogously.

\begin{figure}[h]
\centering
\includegraphics[width=.35\textwidth]{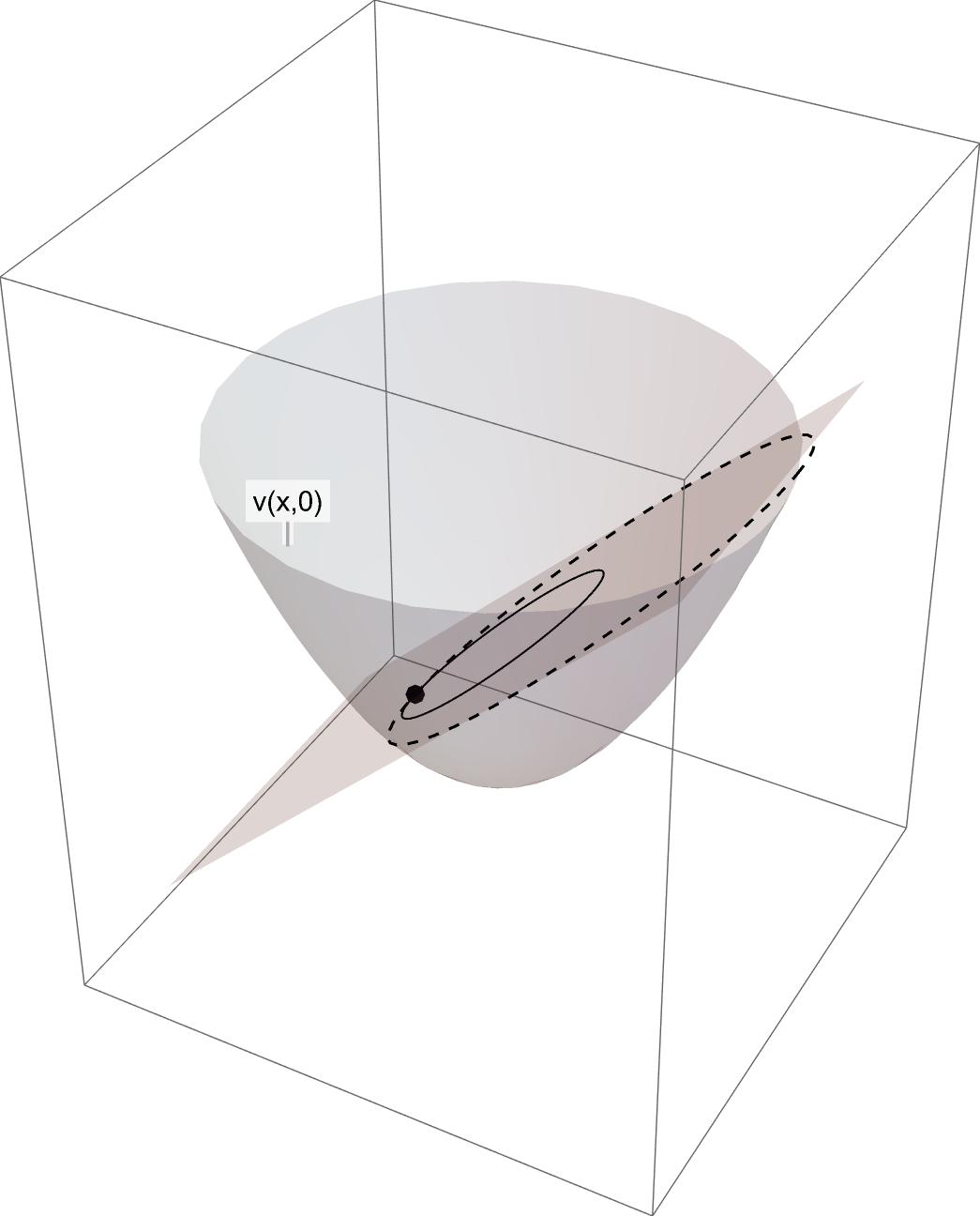}\quad\quad\quad\quad
\includegraphics[width=.35\textwidth]
{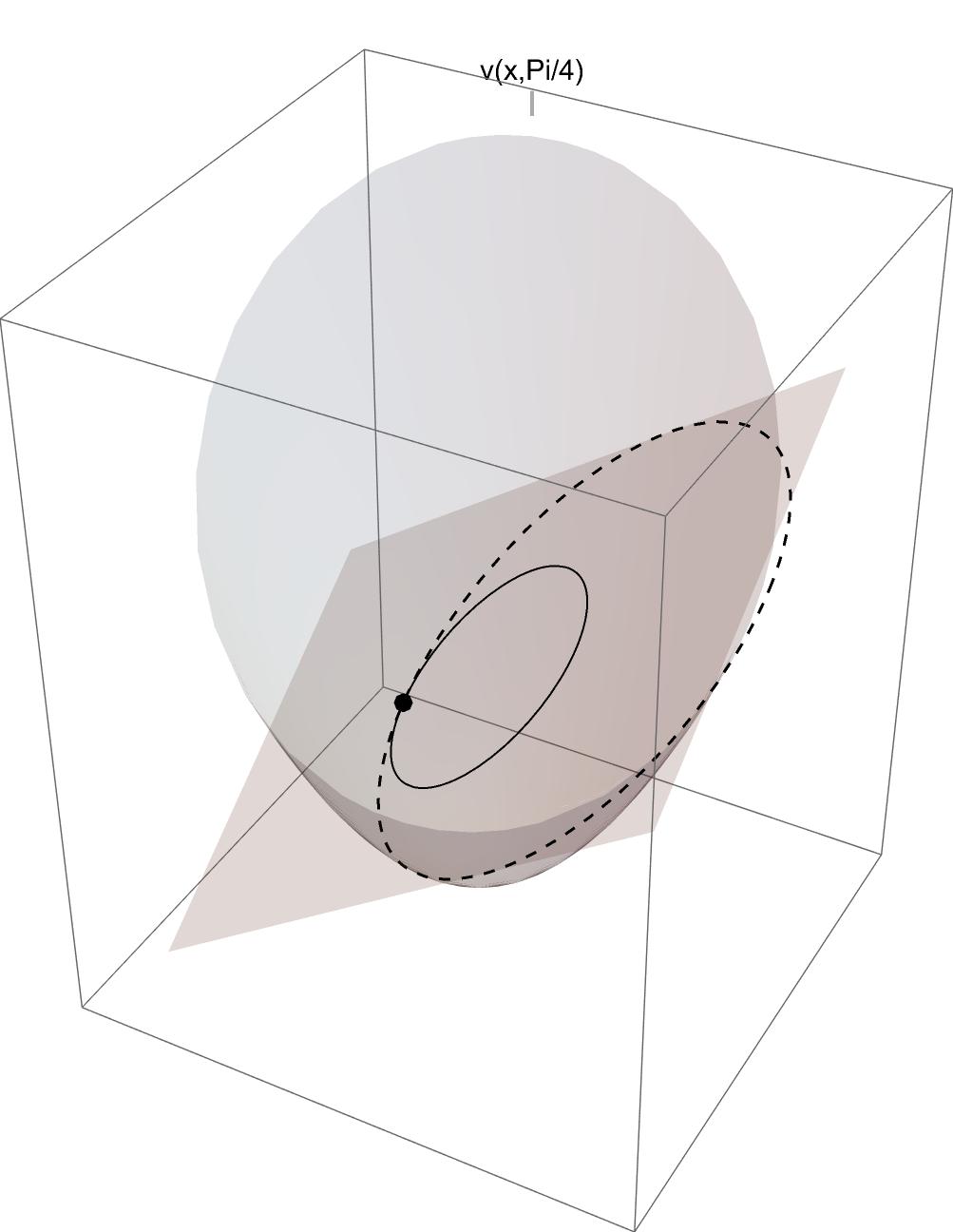}
\caption{}
\label{fig:1}
\end{figure}

\bigskip
 
\noindent
{\bf Acknowledgments.} This research was funded in whole or in part by the Austrian Science Fund (FWF) [10.55776/P36295]. For open access purposes, the author has applied a CC BY public copyright license to any author accepted manuscript version arising from this submission. This reasearch has been partially supported by INdAM(GNAMPA) and by  DOR2340044 UniPd. 

\section{Notations and Preliminaries}
Throughout the paper, if not further specified, $\Omega$ denotes a boun\-ded, open and convex subset of $\R^{n}$, $n\ge 2$ and $\Omega_T:=\Omega\times(0,T)$ is the associated space-time cylinder, $T>0$. 
The parabolic boundary of $\Omega_{T}$ is defined by
\[
    \partial_{\mathcal{P}}\Omega_{T}
    :=
    \left(\overline\Omega\times\{0\}\right)\cup\left(\partial\Omega\times(0,T)\right).
\]
As already mentioned, we frequently use the identification of the class of
Lipschitz continuous functions $C^{0,1}(\overline{\Omega})$ with the space $W^{1,\infty}(\Omega)$. 
For a continuous function $u\colon\overline{\Omega}\to\R$ we recall the definition of the Lipschitz constant 
$$
	[u]_{0,1;\Omega}:=\sup_{x\not=y, x,y\in\Omega}\frac{|u(x)-u(y)|}{|x-y|}<\infty.
$$
Since $\Omega$ is convex, we have $[u]_{0,1;\Omega}=\|\nabla u\|_{L^\infty (\Omega ,\R^n)}$.


\subsection{$R$-uniform convexity and time dependent bounded slope condition}\label{tBSC}

In this subsection we specify the notion of $R$-uniform convexity and the time dependent bounded slope condition which have already been used in the statement of Theorem~\ref{thm:existence-var-sol}.

\begin{definition}\label{def:unif-conv}\upshape
Let $R>0$. A subset $\Omega\subset \R^n$ is called \emph{$R$-uniformly convex}, if for any  $x_o\in\partial \Omega$ there exists a unit vector $\nu_{x_o}\in\R^n$ pointing outside of $\Omega$ such that
\begin{align}\label{eq:unif-conv}
 R\nu_{x_o}\cdot (x-x_o)\le - \tfrac12 |x-x_o|^2
    \qquad\mbox{for any $x\in\partial\Omega$.}   
\end{align}
\end{definition}

\begin{remark}\label{rem:unif-conv}\upshape
The notion of $R$-uniform convexity implies that for any $x\in\partial\Omega$ there holds
\begin{align*}
      |x-(x_o-R\nu_{x_o})|^2
      &=
      |x-x_o|^2 +2R\nu_{x_o}\cdot (x-x_o)+R^2\le R^2.
\end{align*}
For the last inequality we used that $\nu_{x_o}$ points out of $\Omega$ and \eqref{eq:unif-conv}. In particular, we have $\Omega\subset B_R(x_o-R\nu_{x_o})$.
\end{remark}

Since the boundary values $g$  depend on time, we have to introduce a time dependent version of the bounded slope condition. 

\begin{definition}\label{def:bdslope}
\upshape 
A function $g\colon \partial\Omega\times [0,T]\to \R$ satisfies the {\it time dependent bounded slope condition with constant $Q>0$} (hereafter abbreviated as $t-{\rm BSC}_Q$), if for each $x_o\in\partial\Omega$ there exist  functions $w_{x_o}^\pm \colon [0,T]\rightarrow\R^n$  satisfying $\|w_{x_o}^{\pm}\|_{L^\infty ([0,T],\R^n)}\le Q$ such that
\[
    g(x_o,t)+w_{x_o}^-(t)\cdot(x-x_o)
    \le 
    g(x,t)
    \le 
    g(x_o,t)+w_{x_o}^+(t)\cdot(x-x_o)
\]
for any $(x,t)\in\partial\Omega\times[0,T]$.
\end{definition}

\begin{remark}\label{rem:C2}
\upshape In the framework of stationary problems Miranda \cite{Giusti.2003, Miranda:1965} proved that if $\Omega$ is a uniformly convex, bounded domain
and $v\in C^2(\overline{\Omega})$, then
$v|_{\partial\Omega}$ satisfies the bounded slope condition. In the time dependent setting it follows that if $g(x,t)$ is a $C^{2}(\overline{\Omega_T})$ function, then it satisfies the $t-{\rm BSC}_Q$ for some $Q>0$.
\end{remark}

The time dependent bounded slope condition from Definition~\ref{def:bdslope} involves a function $g$ which is defined only on the lateral boundary of the parabolic cylinder. In the next lemma we will see an application for a function $G$ defined on the whole parabolic cylinder and whose boundary values satisfy the $t-{\rm BSC}_Q$.

\begin{lem}\label{lem:tBSC-G}
Let $G\colon \Omega\times [0,T]\to\R$ with 
$\sup_{t\in [0,T]}\Lip G(\cdot,t)\le Q_o$ such that $g:=G\big\vert_{\partial\Omega\times [0,T]}$
satisfies the $t-{\rm BSC}_{Q}$ for some  constant $Q$. Then for every $x_o\in\partial\Omega$ there exist  $\widetilde w_{x_o}^\pm \colon [0,T]\to \R^n$ with $\|\widetilde w_{x_o}^\pm\|_{L^\infty ([0,T],\R^n)}\le Q + Q_o$  such that
for any $(x,t)\in \overline\Omega\times [0,T]$ we have

\begin{align*}
    g(x_o,t)+\widetilde w_{x_o}^{-}(t)\cdot (x-x_o)
    \le
    G(x,t)
    \le
    g(x_o,t)+\widetilde w_{x_o}^{+}(t)\cdot (x-x_o).
\end{align*}
If we additionally assume that $w_{x_o}^\pm\in W^{1,\infty}([0,T],\R^n)$, then
$$
\|(\widetilde w_{x_o}^\pm)'\|_{L^\infty ([0,T],\R^n)}=\|(w_{x_o}^\pm)'\|_{L^\infty ([0,T],\R^n)}.
$$
\end{lem}

\begin{proof}[Proof] Fix a boundary point $x_o\in\partial\Omega$ and denote by $[0,T]\ni t\mapsto w_{x_o}^{+}(t)$ the $\R^n$-valued function from the $t-{\rm BSC}_{Q}$ satisfied by $g$. Recall that $|w_{x_o}^{+}(t)|\le Q_0$ and
\begin{equation*}
    g(x,t)\le g(x_o,t)+w_{x_o}^{+}(t)\cdot (x-x_o)
\end{equation*}
for any $x\in \partial\Omega$ and any $t\in[0,T]$. Next, we fix a time $t\in[0,T]$. 
Without loss of generality we assume that $x_o=0$, that $g(x_o,t)=0$, and that $\Omega\subset \R^{n-1}\times \R_+$.  Consider 
$x=(x',x_n)\in\Omega$. By $(x',y_n)\in \partial\Omega$ we denote the unique point in $\partial\Omega$ with $0\le y_n <x_n$. Then, we obtain
\begin{align*}
    G(x,t)
    &\le
    G(x',x_n,t)-g(x',y_n,t)+ w_{x_o}^+(t)\cdot (x',y_n)\\
    &\le 
    Q_o (x_n-y_n)+w_{x_o}^+(t)\cdot (x',y_n)\\
    &\le
    Q_o x_n+w_{x_o}^+(t)\cdot (x',y_n)\\
    &=\big[\underbrace{Q_oe_n+ w_{x_o}^+(t)}_{:=\widetilde w_{x_o}^+(t)}]\cdot (x',y_n).
\end{align*}
Thus we have shown that $G(x,t)\le \widetilde w_{x_o}^+(t)\cdot x$ for any $x\in\Omega$, that  $g(x_o,t)=0= \big[\widetilde w_{x_o}^+(t)\cdot x\big]\big\vert_{x_o}$, and that $g(x,t)\le  w_{x_o}^+(t)\cdot x
\le \widetilde w_{x_o}^+(t)\cdot x$ for any $x\in\partial\Omega$.  This means that $g$ satisfies the $t-{\rm BSC}_{Q_1}$ with the larger constant $Q_1=Q+Q_o$.
\end{proof}

\subsection{Mollification in time\label{sec:time-mollification}}
In the definition of variational solutions we do not impose any condition on the time derivative. Therefore, in general we cannot insert the variational solution itself as comparison map in the variational
inequality \eqref{eq:variational-inequality}.  Hence, a suitable mollification procedure in time is needed. For a separable
Banach space  $X$, an initial datum $v_{o}\in X$ and an integrability exponent $1\le r\le\infty$, we consider  
$v\in L^{r}(0,T;X)$ and  $h\in(0,T]$ and define the mollification in time by 
\begin{equation}
	[v]_{h}(t):=\mathrm{e}^{-\frac{t}{h}}v_{o}
	+\tfrac{1}{h}\int_{0}^{t}\mathrm{e}^{\frac{s-t}{h}}v(s)\,\mathrm{d}s,\label{eq:time-mollification}
\end{equation}
whenever $t\in[0,T]$.
Later on, we  use for instance $X=L^{r}(\Omega,\R^{N})$ or $X=W^{1,r}(\Omega,\R^{N})$
and the related parabolic spaces $L^{r}(0,T;L^{r}(\Omega,\R^{N}))$, respectively $L^{r}(0,T;W^{1,r}(\Omega,\R^{N}))$.
One particular features of this mollification is that $\left[v\right]_{h}$
solves the ordinary differential equation 
\begin{equation*}
	\partial_{t}[v]_{h}=\tfrac{1}{h}\big(v-[v]_{h}\big)
\end{equation*}
with initial condition $[v]_{h}(0)=v_{o}$. Note that,
since $[v]_{h}$ solves the above ODE, then clearly any common
membership of both $v$ and its regularization $[v]_{h}$
to a Banach space is passed  to the time derivative $\partial_t [v]_{h}$.
The basic properties of the mollification in time are summarized in
the following lemma (cf. \cite[Lemma 2.2]{Kinnunen.2006} and \cite[Appendix B]{Bogelein.2013}
for the proofs).

\begin{lem}\label{lem:moll-props} Suppose $X$ is a separable Banach
space and $v_{o}\in X$. If $v\in L^{r}(0,T;X)$ for some $r\ge1$,
then also $[v]_{h}\in L^{r}(0,T;X)$, and $[v]_{h}\to v$
in $L^{r}\left(0,T;X\right)$ as $h\downarrow0$ if $r<\infty$. Furthermore, for any
$t_{o}\in(0,T]$ there holds 
\[
	\left\Vert [v]_{h}\right\Vert _{L^{r}\left(0,t_{o};X\right)}
	\le
	\left\Vert v\right\Vert _{L^{r}\left(0,t_{o};X\right)}
	+
	\left[\tfrac{h}{r}\left(1-\mathrm{e}^{-\frac{t_{o}r}{h}}\right)\right]^{\frac{1}{r}}
	\left\Vert v_{o}\right\Vert _{X}.
\]
In the case $r=\infty$, the bracket $\left[\ldots\right]^{1/r}$
in the preceding inequality has to be interpreted as $1$. Moreover,
we have $\partial_{t}[v]_{h}\in L^{r}\left(0,T;X\right)$ with
\begin{equation}\label{magic-dt}
	\partial_{t}[v]_{h}=\tfrac{1}{h}\big(v-[v]_{h}\big).
\end{equation}
If in addition $\partial_{t}v\in L^{r}(0,T;X)$ and $v_o=v(0)$, then 
\[
	\partial_{t}[v]_{h}=\tfrac{1}{h}\int_{0}^{t}\mathrm{e}^{\frac{s-t}{h}}\partial_{s}v(s)\ds
\]
and 
\[
	\left\Vert \partial_{t}[v]_{h}\right\Vert _{L^{p}\left(0,T;X\right)}
	\le
	\left\Vert \partial_{t}v\right\Vert _{L^{p}\left(0,T;X\right)}.
\]
Finally, if $v\in C^{0}([0,T];X)$ and $v_o=v(0)$, then also $[v]_{h}\in C^{0}([0,T];X)$,
$[v]_{h}(0)=v_{o}$ and $[v]_{h}\to v$ in $L^{\infty}([0,T];X)$
as $h\downarrow0$. 
\hfill{}$\Box$ \end{lem}

\section{Properties of variational solutions}\label{sec:prop}

In this section we present some properties of variational solutions that are direct consequences of the definition. These properties can be achieved in a fairly general setting, in particular without assuming the bounded slope condition. 

\begin{lem}\label{lem:initial}
Assume that $f\colon\R^n\to \R$ is convex and $g$ satisfies \eqref{ass:g}. Then, any variational solution $u$ in the sense of Definition~\ref{def:var-sol}, resp.~\ref{def:var-sol-rest} fulfills the initial condition $u(0) = g_o$  in the $L^2$-sense, i.e.~we have
$$
	\lim_{h\downarrow 0} \tfrac{1}{h}\int_0^h\|(g-u)(\tau)\|_{L^{2}(\Omega)}^{2}\,\dtau =0.
$$
\end{lem}

\begin{proof}[Proof] In view of assumption~\eqref{ass:g} we find $L>0$ satisfying~\eqref{Def:L}. Since $g$ is an admissible comparison map in~\eqref{eq:variational-inequality}, resp.~\eqref{eq:variational-inequality-rest}, we have
\begin{align*}
	\tfrac{1}{2}\|(g-u)(\tau)\|_{L^{2}(\Omega)}^{2}
	+
	\iint_{\Omega_\tau}f(\nabla u)\,\dx\dt
	& \le
	\iint_{\Omega_\tau}\big[\partial_{t}g(g-u)
	+
	f(\nabla g)\big]\,\dx\dt\\
	&\le
	\iint_{\Omega_\tau}|\partial_{t}g||g-u|\,\dx\dt
	+
	\tau|\Omega| \sup_{|\xi|\le \LL} f(\xi) ,
\end{align*}
for a.e.~$\tau\in [0,T]$. On the left-hand side we discard the non-negative energy term. Then, we  integrate for given $h>0$ with respect to $\tau$ over the interval $[0,h]$. In this way we get
\begin{align*}
	\tfrac{1}{2}\int_0^h\|(g-u)(\tau)\|_{L^{2}(\Omega)}^{2}\,\dtau
	& \le
	h\iint_{\Omega_h}|\partial_{t}g||g-u|\,\dx\dt
	+
	\tfrac12 h^2|\Omega| \sup_{|\xi|\le \LL} f(\xi) .
\end{align*}
We divide both sides by $h$ and apply Young's inequality to the first term on the right. This gives
\begin{align*}
	\tfrac{1}{h}\int_0^h\|(g-u)(\tau)\|_{L^{2}(\Omega)}^{2}\,\dtau
	&\le
	\iint_{\Omega_h} \big[|\partial_{t}g|^2+|g-u|^2\big]\, \dx\dt
	+
	 h|\Omega| \sup_{|\xi|\le \LL} f(\xi) .
\end{align*}
Note that $\partial_tg\in L^2(\Omega_T)$ implies $g\in C^0( [0,T];L^2(\Omega))$.
As the right-hand side vanishes in the limit $h\downarrow 0$ and $g(0)=g_o$, the claim follows.
\end{proof}

In the next lemma we show that variational solutions in the sense of Definition~\ref{def:var-sol} apriori are of class $C^0( [0,T]; L^2(\Omega))$. The  argument is similar to \cite[Lemma 3.2]{BDS}. Here, however,  we have to estimate the difference of energy terms in a different way.

\begin{lem}
Assume that $f\colon\R^n\to \R$ is convex and $g$ satisfies~\eqref{ass:g}. 
Then, every variational solution $u$  in the sense of Definition~\ref{def:var-sol} satisfies 
$$
    u\in C^0\big( [0,T]; L^2(\Omega)\big). 
$$
\end{lem}

\begin{proof}[Proof]
We test the variational inequality for $u$ with 
\begin{equation*}
	v_h= g+[u-g]_h =g+\tfrac1{h}\int_0^t\mathrm e^\frac{s-t}{h}\big(u(s)-g(s)\big)\,\ds,
\end{equation*}
where $[u-g]_h$ is defined as in \eqref{eq:time-mollification} with initial value $v_o=0$. 
In view of Lemma~\ref{lem:moll-props}, one easily shows $v_h\in \Lip_{g}(\Omega_T)$, $\partial_t v_h\in L^2(\Omega_T)$, and $v_h(0)=g_o\in L^2(\Omega)$, so that $v_h$ is admissible in~\eqref{eq:variational-inequality}. In particular, we have a uniform bound for the spatial gradient $|\nabla v_h|$ in $\Omega_T$. In fact, using the definition of $v_h$ and assumption~\eqref{ass:g} we get
\begin{align*}
	|\nabla v_h|
	&= 
    \bigg|\nabla g +\tfrac1{h}\int_0^t\mathrm e^\frac{s-t}{h}\big(\nabla u(s)-\nabla g(s)\big)\,\ds\bigg|\\
	&\le 
    \LL+ \big(\|\nabla u\|_{L^\infty (\Omega_T,\R^n)}+\LL\big) \tfrac1{h}\int_0^t\mathrm e^\frac{s-t}{h}\ds\\
	&\le 
    2\LL + \|\nabla u\|_{L^\infty (\Omega_T,\R^n)}=:\widetilde L.
\end{align*}
Testing the variational inequality~\eqref{eq:variational-inequality} with $v_h$  we obtain
\begin{align*}
	\tfrac{1}{2}\|(v_h-u)(\tau)\|_{L^{2}(\Omega)}^{2}
	& \le
	\iint_{\Omega_\tau}\partial_{t}v_h(v_h-u)\,\dx\dt
	+
	\iint_{\Omega_\tau} \big[f(\nabla v_h)- f(\nabla u)\big]\,\dx\dt
\end{align*}
for a.e.~$\tau\in [0,T]$.  In the last displayed inequality we pass to the supremum over $\tau\in (0,T)$. This implies
\begin{align*}
	\tfrac{1}{2}\sup_{\tau\in(0,T)}& \|(v_h-u)(\tau)\|_{L^{2}(\Omega)}^{2}\\
	& \le
	\sup_{\tau\in(0,T)} \iint_{\Omega_\tau}\partial_{t}v_h(v_h-u)\,\dx\dt
	+
	\iint_{\Omega_T} \big| f(\nabla v_h)- f(\nabla u)\big|\,\dx\dt.
\end{align*}
The second term on the right-hand side can be bounded by means of the gradient bound for $|\nabla v_h|$ in $\Omega_T$. Indeed, we get
\begin{align*}
	\big| f(\nabla v_h)- f(\nabla u)\big|
	&\le
	\sup_{|\xi|\le \widetilde L} |\nabla f(\xi)| | \nabla v_h-\nabla u|
    \to 
    0\quad\mbox{in $L^1(\Omega_T)$ as $h\downarrow 0$,}
\end{align*}
so that
\begin{equation*}
	\lim_{h\downarrow 0}\iint_{\Omega_T} \big| f(\nabla v_h)- f(\nabla u)\big|\,\dx\dt=0.
\end{equation*}
Moreover, exploiting \eqref{magic-dt} we infer
\begin{align*}
    &\iint_{\Omega_\tau}\partial_{t}v_h(v_h-u)\,\dx\dt \\
    &\qquad=
    \iint_{\Omega_\tau}\Big[\partial_{t}g\big(g-u+[u-g]_h\big) +
    \partial_{t}[u-g]_h \big(g-u+[u-g]_h\big)\Big] \,\dx\dt \\
    &\qquad\le
    \iint_{\Omega_\tau} |\partial_{t}g| \big|u-g-[u-g]_h\big| \,\dx\dt .
\end{align*}
so that
\begin{equation*}
	\limsup_{h\downarrow 0} \sup_{\tau\in(0,T)} \iint_{\Omega_\tau}\partial_{t}v_h(v_h-u)\dx\dt\le 0,
\end{equation*}
Inserting the preceding estimates above, we arrive at
\begin{equation*}
	\lim_{h\downarrow 0} \| v_h-u\|_{L^\infty (0,T; L^2(\Omega))}=0.
\end{equation*}
Note that $g\in C^0( [0,T];L^2(\Omega))$, since $\partial_tg\in L^2(\Omega_T)$. This however means that also $v_h$ belongs to
$C^0( [0,T];L^2(\Omega))$. Thus $u\in L^\infty (0,T; L^2(\Omega))$ is the limit with respect to  $L^\infty (
0,T; L^2(\Omega))$-convergence of functions $v_h$ which are themselves of class $C^0( [0,T];L^2(\Omega))$. As a result
the limit map $u$ is continuous with respect to time, i.e.~$u\in C^0( [0,T];L^2(\Omega))$. This proves the claim.
\end{proof}


\section{Existence of solutions for regularized boundary data}
\label{sec:ex-regularized}

In this section we establish for sufficiently regular boundary data $g$ the existence of a variational solution to the gradient constrained problem in the sense of Definition \ref{def:var-sol-rest} via the method of {\it minimizing movements} (finite time discretization). We assume that  $g\in \Lip(\Omega_T)$ additionally satisfies
\begin{equation}\label{g-reg}
	\mbox{$\partial_tg\in L^2(\Omega_T)\cap L^\infty\big( 0,T; W^{1,\infty}(\Omega)\big)$ 
	and $g_o:=g(0)\in W^{1,\infty}(\Omega)$}.
\end{equation}
This implies in particular that $g\in C^0( [0,T], L^2(\Omega))$. Furthermore, let $\LL>0$ such that
\begin{equation}\label{g-constraint}
	\|\nabla g\|_{L^\infty(\Omega_T,\R^n)} < \LL.
\end{equation}

\begin{prop}\label{prop:existence-reg}
Let $f\colon\R^n\to \R$ be convex, $L>0$ and suppose that $g\in \Lip(\Omega_T)$ satisfies \eqref{g-reg} and \eqref{g-constraint}.
Then, there exists a unique variational solution $u$ in $\Omega_T$ in the class $\Lip_g(\Omega_T,\LL)$ in the sense of Definition~\ref{def:var-sol-rest}. Moreover, we have $\partial_t u\in L^2(\Omega_T)$.
\end{prop}

The proof of Proposition~\ref{prop:existence-reg} will be established in  \S\,\ref{sec:ellipticfunc}--\ref{sec:VI-limit}.

\subsection{A sequence of minimizers for variational functionals}\label{sec:ellipticfunc}

We fix a step size
$h\in(0,1]$ and define for $i\in\N_0$ with $ih\le T$ the time-discretized boundary values
$g_i:=g(ih)\in W^{1,\infty} (\Omega)$. Note that $\|\nabla g_i\|_{L^\infty (\Omega,\R^n)}\le L$.
Our goal is to inductively construct a sequence $u_i\in g_i + W^{1,\infty}_0(\Omega)$
of Lipschitz minimizers to certain elliptic variational functionals satisfying the gradient constraint $\|\nabla u_i\|_{L^\infty (\Omega,\R^n)}\le L$. The precise construction is as follows. Suppose that
$u_{i-1}\in g_{i-1}+ W^{1,\infty}_0(\Omega)$ for
some $i\in\N$ has already been defined. If $i=1$, then $u_0=g_o$ is the initial boundary datum. 
Then, we let $u_i$ be the minimizer of the variational functional
\[
	\mathsf F_{i}[v]
	:=
	\int_{\Omega}f(\nabla v)\,\dx +\tfrac{1}{2h}\int_\Omega |v-u_{i-1}|^2\,\dx\, .
\]
in the class of functions $v\in g_i+ W^{1,\infty}_0(\Omega)$ with $\|\nabla v\|_{L^\infty (\Omega,\R^n)}\le L$.
Observe that this class in non-empty since $v=g_i$ is admissible. Note also that $\mathbf F_{i}$ is bounded
on this function class. The existence of a unique minimizer $u_i$ can be
deduced by means of standard compactness arguments (i.e.~by use of the Arzela-Ascoli theorem) using the convexity of $f$.
We note that $\|\nabla u_i\|_{L^\infty (\Omega,\R^n)}\le L$ for any $i\in\N_0$ with $ih\le T$, by construction.

Next, we want to compare the energy of $u_i$ with the one of $u_{i-1}$. This is not directly possible because the boundary values of $u_{i-1}$ do not coincide with those  of $u_i$. We compensate for this by choosing as comparison function $u_{i-1}+g_i-g_{i-1}$ instead of $u_{i-1}$.  From
$$
	\mathsf F_{i}[u_i]\le \mathsf F_{i}[u_{i-1}+g_i-g_{i-1}]
$$
we obtain
\begin{align*}
	\int_{\Omega}f(\nabla u_i)\,\dx &+
	\tfrac{1}{2h}\int_\Omega |u_i-u_{i-1}|^2\,\dx\\
	&\le
	\int_{\Omega}f\big(\nabla (u_{i-1}+g_i-g_{i-1})\big)\,\dx +\tfrac{1}{2h}\int_\Omega |g_i-g_{i-1}|^2\,\dx\\
	&=
	\int_{\Omega}f(\nabla u_{i-1})\,\dx
	+
	\mathbf I_i+\mathbf{II}_i,
 \end{align*}
where
\begin{align*}
	\mathbf I_i&:= \tfrac{1}{2h}\int_\Omega |g_i-g_{i-1}|^2\,\dx,\\
	\mathbf{II}_i
    &:=
    \int_{\Omega}\big[f\big(\nabla (u_{i-1}+g_i-g_{i-1})\big)-f(\nabla u_{i-1})\big]\,\dx.
\end{align*}
For $\mathbf I_i$ we have
\begin{align*}
	\mathbf I_i
	&\le
	\tfrac{1}{2h}\int_\Omega \bigg|\int_{(i-1)h}^{ih} \partial_tg\,\dt\bigg|^2\dx
	\le
	\tfrac12\iint_{\Omega\times [(i-1)h,ih]}|\partial_tg|^2\,\dx\dt.
\end{align*}
Similarly, we obtain
$$
	\big|\nabla (g_i-g_{i-1})\big|
	\le
	\bigg|\int_{(i-1)h}^{ih} \partial_t\nabla g\,\dt\bigg|
	\le
	h\|\partial_t\nabla g\|_{L^\infty (\Omega_T,\R^n)},
$$
so that
$$
	\big\| \nabla (u_{i-1}+g_i-g_{i-1})\big\|_{L^\infty(\Omega)}
	\le
	L+ h\|\partial_t\nabla g\|_{L^\infty (\Omega_T,\R^n)}.
$$
Note that the right-hand side is finite due to assumption~\eqref{g-reg}. 
The preceding inequality can be used to estimate the integral $\mathbf{II}_i$. in fact, we get
\begin{align*}
	|\mathbf{II}_i|
	&\le
	\int_{\Omega}\big|f\big(\nabla (u_{i-1}+g_i-g_{i-1})\big)-f(\nabla u_{i-1})\big|\,\dx\\
	&\le
	\sup_{|\xi|\le L+ \|\partial_t\nabla g\|_{L^\infty (\Omega_T,\R^n)}} 
	|\nabla f(\xi)|
	\int_\Omega |\nabla(g_i-g_{i-1})|\,\dx\\
	&\le
	\mathsf K\iint_{\Omega\times [(i-1)h,ih]}|\partial_t\nabla g|\,\dx\dt,
\end{align*}
where 
\begin{equation}\label{def-K}
    \mathsf K
    :=
    \sup_{|\xi|\le L+ \|\partial_t\nabla g\|_{L^\infty (\Omega_T,\R^n)}} |\nabla f(\xi)|.
\end{equation}
Substituting the estimates for $\mathbf I_i$ and $\mathbf{II}_i$ gives
\begin{align}\label{energy-00}\nonumber
	&\int_{\Omega}f(\nabla u_i)\,\dx +
	\tfrac{1}{2h}\int_\Omega |u_i-u_{i-1}|^2\,\dx\\
	&\qquad\le
	\int_{\Omega}f(\nabla u_{i-1})\,\dx
	+
	\iint_{\Omega\times [(i-1)h,ih]} \big[\tfrac12|\partial_tg|^2 + K |\partial_t\nabla g|\big]\,\dx\dt.
\end{align}
By inductively comparing the energy of $u_i$ with that of $u_{i-1}$ for $i\in \{ 0,\dots ,\ell\}$
in the described way and summing with respect to $i$, we get for any $\ell\in \N$ with $\ell h\le T$ the {\em energy estimate}
\begin{align}\label{energy-0}\nonumber
	&\int_\Omega f(\nabla u_\ell)\,\dx +
	\tfrac1{2h}\sum_{i=1}^\ell\int_\Omega |u_i-u_{i-1}|^2\,\dx\\
	&\qquad\le  
	\int_\Omega f(\nabla g_{o})\,\dx +
    \iint_{\Omega_T}\big[\tfrac12|\partial_tg|^2 + K|\partial_t\nabla g| \big]\,\dx\dt.
\end{align}
The first term of the right-hand side actually contains no new information, because due to the construction of the minimizers $u_\ell$ the $W^{1,\infty}(\Omega) $-norm of $u_\ell$ is bounded anyway. Indeed,
for any $\ell$ as above we have
\begin{align}\label{energy-1}\nonumber
	\|u_\ell \|_{L^\infty (\Omega)}+
	\|\nabla u_\ell\|_{L^\infty (\Omega,\R^n)}
	&\le
	\| g_\ell\|_{L^\infty (\Omega)}+
	\| u_\ell -g_\ell\|_{L^\infty (\Omega)} + 
	\LL \\
	&\le
	 \| g\|_{L^\infty (\Omega_T)}+ 
	 \diam(\Omega)\| \nabla u_\ell -\nabla g_\ell\|_{L^\infty(\Omega)}
	 +\LL\nonumber
	 \\
	 &\le
	 \| g\|_{L^\infty (\Omega_T)}+\LL \big[ 1+2 \diam(\Omega)\big]
	 =:
	 \LL_1.
\end{align}
From now on we consider only such values $h_\ell\in (0,1]$ which satisfy $\ell
:= \frac{T}{h_\ell}\in \N$. Then we define functions $u^{(h_\ell)}\colon\Omega\times (-h_\ell,T]\to
\R$ and $g^{(h_\ell)}\colon\Omega\times (-h_\ell,T]\to
\R$ by
\begin{align*}
	\mbox{$u^{(h_\ell)}(\cdot , t):= u_i$ and $g^{(h_\ell)}(\cdot , t):= g_i$ \quad for $t\in \big((i-1)h_\ell, ih_\ell\big]$ with
	$i\in\{ 0,\dots ,\ell\}$.}
\end{align*}
Note that both $u^{(h_\ell)}$ and $g^{(h_\ell)}$ are  piecewise constant with respect to time and  that
$u^{(h_\ell)}(t)=g^{(h_\ell)}(t)$ on $\partial\Omega$ for any $t\in
(0,T)$ in the usual sense of traces of continuous maps on $\overline\Omega$. This is true since 
both $u^{(h_\ell)}(t)$ and $g^{(h_\ell)}(t) $ can be uniquely extended to
$\overline\Omega$. 
From \eqref{energy-1} we conclude that
\begin{align}\label{energy-2}
	\sup_{t\in [0,T]}\Big[ \big\| u^{(h_\ell)}\big\|_{L^\infty(\Omega)}
	+ \big\| \nabla u^{(h_\ell)}\big\|_{L^\infty(\Omega,\R^n)}\Big]
	&\le 
	\LL_1.
\end{align}
Thus, there exist a map
 $$
 	u\in \bigcap_{q\ge 1}L^{q}(0,T;W^{1,q}(\Omega ))
$$
and a subsequence -- still denoted by  $(u^{(h_\ell)})_{\ell\in\N}$ -- such that
\begin{equation}\label{weak-conv}
\left\{ \begin{array}{cl}
	u^{(h_\ell)}\rightharpoondown u & \mbox{weakly in \ensuremath{L^{q}(\Omega_{T})} for any $q\ge 1$,}\\[5pt]
	\nabla u^{(h_\ell)}\rightharpoondown \nabla u & \mbox{weakly in \ensuremath{L^{q}(\Omega_{T},\R^n)}  for any $q\ge 1$.}
\end{array}\right.
\end{equation}
Observe that by lower-semicontinuity we have for any $q\ge 1$ that
\begin{align}\label{bound-Lq}
	\bigg(\biint_{\Omega_T}|\nabla u|^q\,\dx\dt\bigg)^\frac{1}{q}
	&\le
	\liminf_{\ell\to\infty}
    \bigg(\biint_{\Omega_T}\big|\nabla u^{(h_\ell)}\big|^q\,\dx\dt\bigg)^\frac{1}{q}
	\le\LL,
\end{align}
so that $\| \nabla u\|_{L^\infty(\Omega_T,\R^n)}\le \LL$. In particular, we have $u\in L^\infty( 0,T;W^{1,\infty}(\Omega))$.

So far we have no information about the time derivative of the limit map $u$.
This can be achieved by exploiting the estimate for the second term on the
left-hand side of \eqref{energy-0}. In fact, this term can be translated into an estimate
of the time derivative of the function $\widetilde 
u^{(h_\ell)}\colon \Omega\times(-h_\ell, T]\to\R$, which is constructed by linear interpolation in time  of the minimizing functions $u_{i-1}$ and $u_i$ on the time interval
$( (i-1)h_\ell, ih_\ell]$. 
More precisely,  we define
$$
	\mbox{
	$\widetilde u^{(h_\ell)}(\cdot, t)
	= 
	\big( i-\frac{t}{h_\ell}\big)u_{i-1}
	+ 
	\big(1- i+\frac{t}{h_\ell}\big)u_{i}
	$
	\quad for $t\in\big( (i-1)h_\ell, ih_\ell\big]$
	}
$$
whenever $i\in\{ 0,\dots,\ell\}$. Obviously, $u^{(h_\ell)}$ and
$\widetilde u^{(h_\ell)}$ coincide on $\Omega\times (-h_\ell, 0]$ by definition.
For $t\in ( (i-1)h_\ell, ih_\ell]$ the time derivative of
$\widetilde u^{(h_\ell)}$ can easily be computed as
$$
	\partial_t\widetilde u^{(h_\ell)}(\cdot ,t)=\tfrac{1}{h_\ell} (u_i-u_{i-1}).
$$
We use this  to rewrite the second term on the left-hand side
of \eqref{energy-0}. 
This yields
\begin{align}\label{est-time}
	\tfrac12\iint_{\Omega_T}\big|\partial_t\widetilde u^{(h_\ell)}\big|^2\,\dx\dt\nonumber
	&\le  
	\int_\Omega f(\nabla g_{o})\,\dx
	+
	\iint_{\Omega_T}\big[\tfrac12|\partial_tg|^2 + \mathsf 
    K|\partial_t\nabla g| \big]\,\dx\dt \nonumber \\
	&\le |\Omega|
	\sup_{|\xi|\le\LL} |f(\xi)| + 
	\iint_{\Omega_T}\big[\tfrac12|\partial_tg|^2 + 
    \mathsf K|\partial_t\nabla g| \big]\,\dx\dt.
\end{align} 
Analogous to \eqref{energy-2} we have
\begin{align}\label{energy-3}
	\sup_{t\in [0,T]}\Big[ \big\| \widetilde u^{(h_\ell)}
	\big\|_{L^\infty(\Omega)}
	+ \big\| D\widetilde u^{(h_\ell)}\big\|_{L^\infty(\Omega,\R^n)}\Big]
	&\le \LL_1
	.
\end{align}
Thus, the sequence $\widetilde u^{(h_\ell)}$ is uniformly bounded
in $L^\infty ( 0,T; W^{1,\infty}(\Omega))$. In addition, the sequence of time derivatives 
$\partial_t\widetilde u^{(h_\ell)}$ is  uniformly bounded in $L^2(\Omega_T)$.
Together this implies that the  sequence $(\widetilde u^{(h_\ell)})_{\ell \in \N}$ is
uniformly bounded in $W^{1,2}(\Omega_T)$. By compactness,
there exists a subsequence -- still  denoted
by $(\widetilde u^{(h_\ell)})_{\ell\in\N}$  -- and a limit function $\widetilde u\in \bigcap_{q\ge 1}
L^q(0,T;W^{1,q}(\Omega))$ with $\partial_t\widetilde u\in L^2(\Omega_T)$, such that  in the limit  $\ell\to\infty$ we have
\begin{equation}\label{conv-tilde}
\left\{
\begin{array}{cl}
	\widetilde u^{(h_\ell)}\to \widetilde u &
	\quad \mbox{strongly in $L^{2}(\Omega_T)$,}\\[5pt]
	\widetilde u^{(h_\ell)}\rightharpoondown \widetilde u &
	\quad \mbox{weakly  in $L^q\big(0,T;W^{1,q}(\Omega)\big)$
	for any $q\ge 1$,} \\[5pt]
	\partial_t\widetilde u^{(h_\ell)}\rightharpoondown \partial_t\widetilde u &
	\quad \mbox{weakly in $L^2(\Omega_T)$.}
\end{array}
\right.
\end{equation}
As a first consequence of the weak $L^2$-convergence of the time derivative $\eqref{conv-tilde}_3$ and the lower-semicontinuity
of the $L^2$-norm we conclude from \eqref{est-time} that
\begin{align*}
	\iint_{\Omega_T}|\partial_t\widetilde u|^2\,\dx\dt
	&\le
	\liminf_{\ell\to\infty} 
	\iint_{\Omega_T}\big|\partial_t\widetilde u^{(h_\ell)}\big|^2\,\dx\dt\\
	&\le
	2|\Omega|
	\sup_{|\xi|\le\LL} |f(\xi)| + 
	\iint_{\Omega_T}\big[|\partial_tg|^2 + 2\mathsf K|\partial_t\nabla g| \big]\,\dx\dt.
\end{align*}
In the further course of the proof we will see that $\widetilde u=u$, so that  $u$ has a time derivative in $L^2(\Omega_T)$
with the quantitative estimate
\begin{align}\label{est-time-u}
	\iint_{\Omega_T}|\partial_t u|^2 \,\dx\dt
	&\le
	2|\Omega|
	\sup_{|\xi|\le\LL} |f(\xi)| + 
	\iint_{\Omega_T}\big[|\partial_tg|^2 + 2\mathsf K|\partial_t\nabla g| \big]\,\dx\dt.
\end{align}
Now we establish that $\widetilde u= u$. 
By comparing the two sequences $u^{(h_\ell)}$ and $\widetilde u^{(h_\ell)}$, it is easy to see that their limits coincide.  Indeed, we have
$$
	\mbox{$\big|\widetilde u^{(h_\ell)}-u^{(h_\ell)}\big|
	\le
	|u_i-u_{i-1}|$
	\quad for $t\in \big((i-1)h_\ell, ih_\ell\big]$,}
$$
which together with \eqref{energy-0} leads us to
\begin{align*}
	\iint_{\Omega_T}\big|\widetilde u^{(h_\ell)}-u^{(h_\ell)}\big|^2\,\dx\dt
	&\le
	h_\ell\sum_{i=1}^\ell\int_\Omega|u_i-u_{i-1}|^2\,\dx\\
	&\le 
	h_\ell^2\bigg[
	2|\Omega|
	\sup_{|\xi|\le\LL} |f(\xi)| + 
	\iint_{\Omega_T}\big[|\partial_tg|^2 + 2 \mathsf K|\partial_t\nabla g| \big]\,\dx\dt\bigg].
\end{align*}
Taking into account \eqref{conv-tilde}$_1$, we obtain  that
$u^{(h_\ell)}\to\widetilde u$ strongly in $L^2(\Omega_T)$ as $\ell\to\infty$.
Therefore, we have $\widetilde u=u$. The strong $L^2(\Omega_T)$-convergence allows us to pass to another (not relabelled) subsequence, which then converges almost everywhere, i.e.~$u^{(h_\ell)}\to u$ a.e.~in
$\Omega_T$ as $\ell\to\infty$. 

Now we turn our attention to the initial time $t=0$. With 
\eqref{conv-tilde}$_1$ and \eqref{est-time} we obtain for $t\in(0,T)$ that
\begin{align*}
	\frac1t\iint_{\Omega\times [0,t]}&|u(\tau)-g_o|^2\,\dx\dtau\\
	&=
	\lim_{\ell\to\infty}\frac1t\iint_{\Omega\times [0,t]}\big|\widetilde u^{(h_\ell)}(\tau)-g_o\big|^2\,\dx\dtau \\
	&=
	\lim_{\ell\to\infty}\frac1t\iint_{\Omega\times [0,t]} 
	\bigg|\int_0^\tau\partial_\tau\widetilde u^{(h_\ell)}(s)\,\ds\bigg|^2\,\dx\dtau \\
	&\le
	\lim_{\ell\to\infty}\frac1t\iint_{\Omega\times [0,t]} 
	\tau\int_0^\tau\big|\partial_\tau\widetilde u^{(h_\ell)}\big|^2\,\ds\dx\dtau \\
	&\le
	\bigg[ 
	2|\Omega|
	\sup_{|\xi|\le\LL} |f(\xi)| + 
	\iint_{\Omega_T}\big[|\partial_tg|^2 + 2\mathsf K|\partial_t\nabla g| \big]\,\dx\dt\bigg] 
	\frac1t \int_0^t\tau\,\dtau\\
	&= 
	t \bigg[ 
	|\Omega|
	\sup_{|\xi|\le\LL} |f(\xi)| + 
	\iint_{\Omega_T}\big[\tfrac12|\partial_tg|^2 + \mathsf K|\partial_t\nabla g| \big]\,\dx\dt\bigg].
\end{align*}
This implies
\begin{align*}
	\lim_{t\downarrow 0}\frac1t\iint_{\Omega\times [0,t]}|u(\tau)-g_o|^2\,\dx\dtau
	=
	0.
\end{align*}

\subsection{Minimizing properties of the approximation sequences}\label{sec:min-prop} In this subsection we show that the piecewise constant functions $u^{(h_\ell)}$ constructed via time discretization minimize a certain integral functional on a space-time cylinder. In fact, for fixed $\tau\in[0,T]$ and $\ell\in\N$ the function $u^{(h_\ell)}$ minimizes 
$$
 	\mathsf F^{(h_\ell)}[v]
	:=
	\iint_{\Omega_\tau}
	\Big[f\big(\nabla u^{(h_\ell)}\big) +
    \tfrac{1}{2h_\ell}\big|v(t)-u^{(h_\ell)}(t-h_\ell)\big|^2 \Big]\,\dx\dt 
$$
in the class of functions
\begin{equation}\label{comp-maps}
	\mbox{$v\in g^{(h_\ell)}+L^\infty \big(0,\tau;W^{1,\infty}_{0}(\Omega)\big)$
	with $\displaystyle \| Dv\|_{L^\infty (\Omega_\tau)}
	\le L$.}
\end{equation}
This results from a simple calculation, using the definition of  $u^{(h_\ell)}$, the minimizing property of $u_i$, and the  definition of the functional
$\mathbf F^{(h_\ell)}$. Indeed, denoting by $\boldsymbol \chi_{[0,\tau]}$ the characteristic function of the interval $[0,\tau]$, we have
\begin{align*}
	\mathsf F^{(h_\ell)}\big[u^{(h_\ell)}\big]
	&=
	\sum_{i=1}^{\ell} \int_{(i-1)h_\ell}^{ih_\ell} \boldsymbol \chi_{[0,\tau]}
	\int_{\Omega} \Big[\tfrac{1}{2h_\ell}|u_i-u_{i-1}|^2 + f(\nabla u_i)\Big]
	\,\dx \dt\\
	&=
	\sum_{i=1}^{\ell} \int_{(i-1)h_\ell}^{ih_\ell} \boldsymbol \chi_{[0,\tau]}\mathsf F_i[u_i] \,\dt 
	\le
	\sum_{i=1}^{\ell} \int_{(i-1)h_\ell}^{ih_\ell}  \boldsymbol\chi_{[0,\tau]}\mathsf F_i[v(t)] \,\dt \\
	&=	
	\sum_{i=1}^{\ell} \int_{(i-1)h_\ell}^{ih_\ell} \chi_{[0,\tau]}
	\int_{\Omega}
	\Big[
    f(\nabla v(t))
    +
    \tfrac{1}{2h_\ell}
	\big|v(t)-u^{(h_\ell)}(t-h_\ell)\big|^2\Big]\,\dx\dt \\
	&=	
	\mathsf F^{(h_\ell)}[v].
\end{align*}
We use this minimality condition below to derive a variational inequality for the approximating functions $u^{(h_\ell)}$.  In this discrete variational inequality we can then proceed to the limit $\ell\to\infty$ and obtain that $u$ is the desired solution. The precise procedure is as follows. We first re-write the above inequality in terms of the total energy of $u^{(h_\ell)}$ as follows
\begin{align*}
	\iint_{\Omega_\tau} &f\big(\nabla u^{(h_\ell)}\big) \,\dx\dt\\
	&\le 
	\iint_{\Omega_\tau}f(\nabla v) \,\dx\dt \\
	&\phantom{\le\,}+
	\tfrac{1}{2h_\ell} \iint_{\Omega_\tau}
	\Big[\big|v(t)-u^{(h_\ell)}(t{-}h_\ell)\big|^2 -\big |u^{(h_\ell)}(t)-
	u^{(h_\ell)}(t{-}h_\ell)\big|^2\Big]\,\dx\dt \\
	&= 
	\iint_{\Omega_\tau}f(\nabla v) \,\dx\dt\\ 
	&\phantom{\le\,}+
	\tfrac{1}{h_\ell} \iint_{\Omega_\tau}
	\Big[\tfrac12\big|v-u^{(h_\ell)}\big|^2 + \big(v-u^{(h_\ell)}\big)\big(u^{(h_\ell)}-u^{(h_\ell)}(t{-}h_\ell)\big)\Big]\,\dx\dt. 
\end{align*}
We now replace $v$ by the convex combination of $u^{(h_\ell)}$ and a general admissible function $v$ as in \eqref{comp-maps}, i.e.~the function $w^{(h_\ell)}:=u^{(h_\ell)}+s(v-u^{(h_\ell)})$ with $s\in (0,1)$. Note that $w^{(h_\ell)}$ is still admissible, because
$w^{(h_\ell)}$ fulfills both requirements in \eqref{comp-maps}. 
Using also the convexity of $f$, we obtain
\begin{align*}
	&\iint_{\Omega_\tau} f\big(\nabla u^{(h_\ell)} \big) \,\dx\dt \\
	&\qquad\le 
	\iint_{\Omega_\tau} f\big(\nabla u^{(h_\ell)}+s(\nabla v-\nabla u^{(h_\ell)}) \big) \,\dx\dt\\
	&\qquad\phantom{=\ }+
	\tfrac{1}{h_\ell} \iint_{\Omega_T}
	\Big[\tfrac{s^2}{2}\big|v-u^{(h_\ell)}\big|^2
	+ s\big(v-u^{(h_\ell)}\big)\big(u^{(h_\ell)}-u^{(h_\ell)}(t-h_\ell)\big)\Big]	\,\dx\dt \\
	&\qquad\le 
	\iint_{\Omega_\tau} \Big[(1-s)f\big(\nabla u^{(h_\ell)}\big) + sf(\nabla v)\Big] \,\dt\\
	&\qquad\phantom{=\ }+
	\tfrac{1}{h_\ell}\iint_{\Omega_\tau}
	\Big[\tfrac{s^2}{2}\big|v-u^{(h_\ell)}\big|^2
	+ s\big(v-u^{(h_\ell)}\big)\big(u^{(h_\ell)}-u^{(h_\ell)}(t-h_\ell)\big)\Big]\,\dx\dt.
\end{align*}
Here we absorb the first integral on the right-hand side into the left. Then we divide by $s>0$ and let $s\downarrow 0$.
This yields
\begin{align*}
	\iint_{\Omega_\tau} f\big(\nabla u^{(h_\ell)} \big) \,\dx\dt
	&\le 
	\iint_{\Omega_\tau} f(\nabla v)\,\dx\dt
	+
	\iint_{\Omega_\tau}
	\big(v-u^{(h_\ell)}\big)\Delta_{-h_\ell} u^{(h_\ell)}\,\dx\dt,
\end{align*}
where 
$$
    \Delta_{-h_\ell} u^{(h_\ell)}
    :=
    \tfrac{1}{h_\ell}\big(u^{(h_\ell)}-u^{(h_\ell)}(t-h_\ell)\big).
$$
In the second integral on the left-hand side we perform a discrete integration by parts.
For this, however, it is necessary that $v$ is also defined for negative times. Therefore, we define $v(t):=v(0)$
for $t<0$ and assume  that $v(0)\in L^2(\Omega)$ exists.
This allows to conclude 
\begin{align}\label{approx-minimizer}\nonumber
	\iint_{\Omega_\tau} & f\big(\nabla u^{(h_\ell)} \big)\,\dx\dt \\\nonumber
	&\le 
	\iint_{\Omega_\tau}f(\nabla v) \,\dx\dt +
	\tfrac{1}{h_\ell} \iint_{\Omega_\tau}
	\big(v-u^{(h_\ell)}\big)\big(v-v(t-h_\ell)\big) \,\dx\dt \\\nonumber
	&\quad+
	\tfrac{1}{2h_\ell} \iint_{\Omega_\tau}
	\Big[\big|v-u^{(h_\ell)}\big|^2(t-h_\ell) - 
	\big|v-u^{(h_\ell)}\big|^2\Big] \,\dx\dt \\\nonumber
	&\quad-
	\tfrac{1}{2h_\ell} \iint_{\Omega_\tau}
	\big|v - v(t-h_\ell) - u^{(h_\ell)} + u^{(h_\ell)}(t-h_\ell)\big|^2
        \,\dx\dt \\\nonumber
	&\le
	\iint_{\Omega_\tau}f(\nabla v) \,\dx\dt +
	\iint_{\Omega_\tau}
	\big(v-u^{(h_\ell)}\big)\Delta_{-h_\ell}v \,\dx\dt \\
	&\phantom{=\ }-
	\tfrac{1}{2h_\ell}  \iint_{\Omega\times [\tau-h_\ell,\tau]}
	\big|v-u^{(h_\ell)}\big|^2 \,\dx\dt +
    \tfrac{1}{2h_\ell}  \iint_{\Omega\times [-h_\ell,0]}
	|v-g_o|^2 \,\dx\dt\,.
\end{align}

\subsection{Variational inequality for the limit map}\label{sec:VI-limit}
Our aim here is to perform the limit $h_\ell\downarrow 0$ in \eqref{approx-minimizer}. 
To this aim we need to replace the
boundary condition $v = g^{(h_\ell)}$ on the lateral boundary by the $h_\ell$-independent condition
$v = g$. Therefore we consider $v\in \Lip_g(\Omega_T,\LL)$ satisfying $\partial_tv\in L^2(\Omega_T)$ and $v(0)\in L^2(
\Omega)$. Moreover, we extend $v$ to negative times $t<0$ by  
$v(t):=v( 0)\in L^2(\Omega)$. In \eqref{approx-minimizer} we would like to insert $v+g^{(h_\ell)}-g$ as comparison map. However, this is not possible due to the gradient constraint.
This problem can be compensated by considering a convex combination of $v$ and $g$ instead of $v$. More precisely, for 
$\lambda\in (0,1)$ we consider
$$
	v^{(h_\ell)} := \lambda g +(1-\lambda)v+g^{(h_\ell)}-g
$$
as comparison map. Note that $v^{(h_\ell)}\in g^{(h_\ell)} +L^\infty(0,T;W^{1,\infty}_0(\Omega))$. Moreover, in view of \eqref{g-constraint} we have
\begin{align*}
	\|\nabla (\lambda g +(1-\lambda)v)\|_{L^\infty(\Omega_T,\R^n)}
	\le
	\lambda \|\nabla g\|_{L^\infty(\Omega_T,\R^n)} + 
	(1-\lambda) \|\nabla v\|_{L^\infty(\Omega_T,\R^n)} 
	<
	L, 
\end{align*}
and also
\begin{align*}
	\big\| \nabla g^{(h_\ell)}-\nabla g\big\|_{L^\infty(\Omega_T,\R^n)}
	&\le h_\ell\| \partial_t\nabla g\|_{L^\infty(\Omega_T,\R^n)}.
\end{align*}
Therefore, for $\ell\in\N$ large enough, i.e.~for $\ell$ with 
$$
	h_\ell\| \partial_t\nabla g\|_{L^\infty(\Omega_T,\R^n)}
	\le 
	L- 
	\|\nabla \big(\lambda g+ (1-\lambda)v\big)\|_{L^\infty(\Omega_T,\R^n)},
$$
the function $v^{(h_\ell)}$ fulfills the gradient constraint, so that $v^{(h_\ell)}\in \Lip_g(\Omega_T,\LL)$ and hence it is admissible in \eqref{approx-minimizer}.
In the sequel we consider the individual terms on the right-hand side of \eqref{approx-minimizer}. We start with
the integral involving the integrand $f$. The comparison of the energies of $v^{(h_\ell)}$  and $\lambda g +(1-\lambda)v$ shows
\begin{align*}
	\iint_{\Omega_\tau} &\big| f\big(\nabla v^{(h_\ell)}\big)- 
	f\big( \nabla (\lambda g +(1-\lambda)v)\big)\big|\,\dx\dt\\
	&\le 
    \sup_{|\xi|\le L+ \| \partial_t\nabla g\|_{L^\infty(\Omega_T,\R^n)}}|\nabla f(\xi)|
	\iint_{\Omega_\tau}\big| \nabla g^{(h_\ell)}-\nabla g\big|\,\dx\dt\\
	&\le \mathsf
	K h_\ell \iint_{\Omega_T}| \partial_t \nabla  g|\,\dx\dt , 
\end{align*}
for any $\tau\in[0,T]$, where $\mathsf K$ is defined in \eqref{def-K}. 
From this inequality and the convexity of $f$ we obtain
\begin{align}\label{est:comp-en}
	\iint_{\Omega_\tau}& f\big(\nabla v^{(h_\ell)}\big)\,\dx\dt\nonumber\\
	&\le 
	\iint_{\Omega_\tau} f\big( \nabla (\lambda g +(1-\lambda)v)\big)\,\dx\dt+
	\mathsf K h_\ell \iint_{\Omega_T}| \partial_t \nabla  g|\,\dx\dt \nonumber\\
	&\le
	\lambda \iint_{\Omega_\tau} f(\nabla g)\,\dx\dt +
	(1-\lambda) \iint_{\Omega_\tau} f(\nabla v)\,\dx\dt +
	\mathsf K h_\ell \iint_{\Omega_T}| \partial_t \nabla  g|\,\dx\dt.
\end{align}
Next we consider the term involving the time derivative, i.e~the second integral on the right-hand side of~(\ref{approx-minimizer}), and 
observe that 
\begin{align*}
  \Delta_{-h_\ell} v^{(h_\ell)}
  \to
  \lambda\partial_t g +(1-\lambda)\partial_t v
  \quad\mbox{strongly in }L^2(\Omega_T),
\end{align*}
since $\partial_tv, \partial_tg\in L^2(\Omega_T)$ by assumption. 
Together with~(\ref{weak-conv})$_1$ (here we only need the conclusion for $q=2$),
this implies 
\begin{align}\nonumber \label{time-term-converge}
    \lim_{\ell\to\infty}
    \iint_{\Omega_\tau}&
	\big(v^{(h_\ell)}-u^{(h_\ell)}\big)\Delta_{-h_\ell}v^{(h_\ell)} \,\dx\dt\\
    &=
    \lambda \iint_{\Omega_\tau}(v-u) \partial_t g \,\dx\dt
    +
    (1-\lambda)\iint_{\Omega_\tau}(v-u) \partial_tv\,\dx\dt.
\end{align}
Next, we turn our attention to the last two integrals on
the right-hand side of \eqref{approx-minimizer}. Using  $v^{(h_\ell)}(t)=\lambda g_o+(1-\lambda) v( 0)$
for $t\in(-h_\ell,0)$, the last integral in~(\ref{approx-minimizer}) takes the form
\begin{align}\label{initial-converge}
	\tfrac{1}{2h_\ell}  \iint_{\Omega\times [-h_\ell,0]}\big|v^{(h_\ell)}-u_o\big|^2 \,\dx\dt
	&=
	\tfrac{1}{2}\int_{\Omega}\big|\lambda g_o+(1-\lambda)v(0)-g_o\big|^2 \,\dx\nonumber\\
	&=
	\tfrac{1}{2}(1-\lambda)^2\int_{\Omega} |v(0)-g_o|^2 \,\dx.
\end{align}
It remains to consider the second last term on the right-hand side of ~(\ref{approx-minimizer}), i.e.~the integral
$$
	-\tfrac{1}{2h_\ell}\iint_{\Omega\times[\tau-h_\ell, \tau]} \big| v^{(h_\ell)} - u^{(h_\ell)} \big|^2\,\dx\dt.
$$
Since $\partial_t g\in L^2(\Omega_T)$, we have 
\begin{align*}
	\lim_{\ell\to\infty}
  	\tfrac{1}{2h_\ell} \iint_{\Omega\times [\tau-h_\ell,\tau ]}|g^{(h_\ell)}-g|^2 \,\dx\dt
    \le 
    \lim_{\ell\to\infty}
  	\tfrac{1}{2} \iint_{\Omega\times [\tau-h_\ell,\tau ]}|\partial_t g|^2 \,\dx\dt
  	=0.
\end{align*}
Moreover, since $g$ and $v$ are of class $C^0([0,T];L^2(\Omega))$, it follows that
\begin{align*}
	\lim_{\ell\to\infty}
  	\tfrac{1}{2h_\ell} \iint_{\Omega\times [\tau-h_\ell,\tau ]}|g-g(\tau)|^2 +|v-v(\tau)|^2 \,\dx\dt
  	=0.
\end{align*}
Finally, in view of \eqref{energy-00} and the construction of $\widetilde u^{(h_\ell)}$, we find
\begin{align*}
	\lim_{\ell\to\infty} &
  	\tfrac{1}{2h_\ell} \iint_{\Omega\times [\tau-h_\ell,\tau ]}|u^{(h_\ell)}-\widetilde u^{(h_\ell)}(\tau)|^2 \,\dx\dt \\
    &\le 
    \lim_{\ell\to\infty}
  	h_\ell\bigg[\int_{\Omega\times\{\tau-h_\ell\}}f(\nabla u^{(h_\ell)})\,\dx +
	\iint_{\Omega\times [\tau-2h_\ell,\tau]} \big[\tfrac12|\partial_tg|^2 + \mathsf K |\partial_t\nabla g|\big]\,\dx\dt \bigg] \\
    &\le
    \lim_{\ell\to\infty}
  	h_\ell\bigg[|\Omega|\sup_{|\xi|\le L}|\nabla f(\xi)| +
	\iint_{\Omega_T} \big[\tfrac12|\partial_tg|^2 +\mathsf K |\partial_t\nabla g|\big]\,\dx\dt \bigg]
  	=0.
\end{align*}
This allows us to replace $g$ and $v$ by their restriction to the time slice $\Omega\times\{ \tau\}$, i.e.~by 
$g(\tau)$ and $v(\tau)$, and $u^{(h_\ell)}$ by $u^{(h_\ell)}(\tau)$. It therefore remains
to treat the integral
$$
	 \int_{\Omega}\big|\lambda g(\tau)+(1-\lambda)v(\tau)-\widetilde u^{(h_\ell)}(\tau)\big|^2 \,\dx
$$
in the limit $h_\ell\downarrow 0$. For this aim we need to identify the limit of  $\widetilde u^{(h_\ell)}(\tau)$. We claim that $\widetilde u^{(h_\ell)}(\tau)\rightharpoondown u(\tau)$ weakly in $L^2(\Omega)$. Indeed, observing that
$\widetilde u^{(h_\ell)}(0)=g_o$, we conclude for any $\eta\in L^2(\Omega)$ that
\begin{align*}
	\lim_{\ell\to\infty}
	\int_{\Omega}\widetilde u^{(h_\ell)}(\tau)\eta\,\dx 
    &=
    \lim_{\ell\to\infty}
	\iint_{\Omega_\tau}\partial_t \widetilde u^{(h_\ell)}\eta\,\dx\dt +
	\int_{\Omega}g_o\eta\,\dx\\
    &=
    \iint_{\Omega_\tau}\partial_t u\,\eta\,\dx\dt +
	\int_{\Omega}g_o\eta\,\dx\\
    &=
    \int_{\Omega}u(\tau)\eta\, \dx.
\end{align*}
By lower semicontinuity
we therefore have
\begin{align}\label{final-values-converge}
	\int_{\Omega}\big|\lambda g(\tau)&+(1-\lambda)v(\tau)-u(\tau)\big|^2 \dx\nonumber\\
	&\le
	\liminf_{\ell\to\infty}
	 \int_{\Omega}\big|\lambda g(\tau)+(1-\lambda)v(\tau)-\widetilde u^{(h_\ell)}(\tau)\big|^2 \dx.
\end{align} 
Now we put everything together, i.e.~we first use  \eqref{est:comp-en}, (\ref{time-term-converge}), \eqref{initial-converge}, and (\ref{final-values-converge}) in \eqref{approx-minimizer} to control the terms on the right-hand side,
and then let $\ell\to\infty$. In this last argument, we use the fact that the integral functional related to the convex function $f$  is lower semicontinuous w.r.t.~the convergence $\nabla u^{(h_\ell)}\rightharpoondown \nabla u$
weakly in $L^{q}(\Omega_{T},\R^n)$ for any $q\ge 1$. With these arguments we obtain
\begin{align*}
	\iint_{\Omega_\tau} &f(\nabla u)\,\dx\dt \\
	&\le
	(1-\lambda)\bigg[ \iint_{\Omega_\tau} f(\nabla v) \,\dx\dt +
	\iint_{\Omega_\tau} \partial_t v(v-u) \,\dx\dt\bigg] \\
	&\phantom{\le \ }
	+\lambda\bigg[ \iint_{\Omega_\tau} f(\nabla g) \,\dx\dt +
	\iint_{\Omega_\tau} \partial_t g(v-u) \,\dx\dt\bigg] \\
	&\phantom{=\ }-
	\tfrac12 \big\|(\lambda g +(1-\lambda)v-u)(\tau)\big\|^2_{L^2(\Omega)} 
	+
	\tfrac12(1-\lambda)^2 \|v(0)-g_o\|^2_{L^2(\Omega)} . 
\end{align*}
Letting $\lambda\downarrow 0$ gives
\begin{align*}
	\iint_{\Omega_\tau} f(\nabla u)\,\dx\dt 
	&\le
	\iint_{\Omega_\tau} f(\nabla v) \,\dx\dt +
	\iint_{\Omega_\tau} \partial_t v(v-u) \,\dx\dt\\
	&\phantom{=\ }-
	\tfrac12 \|(v-u)(\tau)\|^2_{L^2(\Omega)} 
	+
	\tfrac12\|v(0)-g_o\|^2_{L^2(\Omega)} . 
\end{align*}
This inequality applies to all
$v\in g+ L^\infty(0,T;W^{1,\infty}_{0}(\Omega))$
satisfying the requirements 
$$
    \| \nabla v\|_{L^\infty (\Omega_T ,\R^n)}
    \le 
    L,\quad
    \partial_t v\in L^2(\Omega_T), 
    \quad\mbox{and  $v(0)\in L^2(\Omega)$. }
$$
Hence, $u$ is a variational solution of the gradient constrained problem. The uniqueness can be deduced as in~\cite[Lemma~3.3]{BDS}. This finishes the proof of Proposition~\ref{prop:existence-reg}.

\section{Existence of solutions to the gradient constrained problem}\label{sec:constr-data}
In this section we prove the existence of a variational solution to the gradient constrained problem in the sense of Definition \ref{def:var-sol-rest} under the assumptions of Theorem~\ref{thm:existence-var-sol}. In particular, we assume the boundary data $g$ to admit merely the regularity given in hypothesis~\eqref{ass:g}. 

\begin{prop}\label{prop:ex-constr}
Let $f\colon\R^n\to \R$ be convex, $L>0$ and suppose that $g$ satisfies \eqref{ass:g} and \eqref{Def:L}.
Then, there exists a unique variational solution $u$ in $\Omega_T$ in the class $\Lip_g(\Omega_T,\LL)$ in the sense of Definition~\ref{def:var-sol-rest}. 
\end{prop}

For the proof of Proposition~\ref{prop:ex-constr} we will rely on the existence result from Proposition~\ref{prop:existence-reg} for more regular boundary data. Therefore, we first regularize the boundary data in such a way that Proposition~\ref{prop:existence-reg} is applicable. Subsequently we will show that the limit function is a solution to the original problem.

\subsection{Construction of regularized data}  
We define {\em regularized boundary values} $g_i$, $i\in\N$, according to \eqref{eq:time-mollification} with
$(g_o, g,h_i)$ instead of $(v_o, v,h)$. Recall that
\begin{equation}\label{def:gi}
	g_i(t):= [g]_{h_i}(t)= \mathrm{e}^{-\frac{t}{h_i}}g_{o}+
    \tfrac{1}{h_i}\int_0^t\mathrm{e}^\frac{s-t}{h_i}g(s)\,\ds.
\end{equation}
Since $g\in L^q\big( 0,T;W^{1,q}(\Omega)\big)$ and $g_o\in W^{1,q}(\Omega)$ for any  $q\ge 1$,
we may use the assertions of Lemma \ref{lem:moll-props} with $r=q$ and $X=W^{1,q}(\Omega)$. From 
Lemma \ref{lem:moll-props} we get
$g_i\in L^{q}(0,T;W^{1,q}(\Omega))$
with $\partial_{t}g_i\in L^{2}(\Omega_{T})$.
 Moreover, we have
\begin{equation}\label{reg-gi}
	\mbox{$\partial_{t}g_i=\tfrac{1}{h_i}(g-g_i)\in L^2(\Omega_T)$
	and
	$\partial_t\nabla g_i=\tfrac{1}{h_i}(\nabla g-\nabla g_i)\in L^q(\Omega_T,\R^n)$.}
\end{equation}
The second assertion follows, as the right-hand side of  the identity for $\partial_{t}g_i$ belongs to $L^{q}\big(0,T;W^{1,q}(\Omega)\big)$.
Since $g_i(0)=g_{o}$, we also have 
$\partial_{t}g_i(0)=\frac{1}{h_i}(g_{o}-g_i(0))=0$.
From the second part of Lemma \ref{lem:moll-props} we obtain
\begin{equation}\label{dtgi}
    \|\partial_t g_i\|_{L^2(\Omega_{t_o})}
    \le 
    \|\partial_t g\|_{L^2(\Omega_{t_o})}
    \quad\mbox{for any $t_o\in(0,T]$.}
\end{equation}
Further, using the convexity of $\xi\mapsto |\xi|^q$ and Jensen's inequality we conclude  from \eqref{def:gi} that
\begin{align*}
	|\nabla g_i(t)|^q
	&=
	\bigg| 
	\mathrm{e}^{-\frac{t}{h_i}}\nabla g_{o}
	+\frac{1-\mathrm{e}^{-\frac{t}{h_i}}}{h_i(1-\mathrm{e}^{-\frac{t}{h_i}})}
	\int_0^t\mathrm{e}^\frac{s-t}{h_i}\nabla g(s)\,\ds
	\bigg|^q\\
	&\le 
	\mathrm{e}^{-\frac{t}{h_i}}|\nabla g_o|^q
	+ 
	\big( 1-\mathrm{e}^{-\frac{t}{h_i}}\big)
	\bigg|\tfrac{1}{h_i(1-\mathrm{e}^{-\frac{t}{h_i}})}
	\int_0^t\mathrm{e}^\frac{s-t}{h_i}\nabla g(s)\,\ds\bigg|^q\\
	&\le 
	\mathrm{e}^{-\frac{t}{h_i}}|\nabla g_o|^q
	+ 
	\tfrac{1}{h_i}
	\int_0^t\mathrm{e}^\frac{s-t}{h_i}|\nabla g(s)|^q\,\ds\\
	&\le 
	\bigg[ \mathrm{e}^{-\frac{t}{h_i}}\|\nabla g_o\|^q_{L^\infty(\Omega,\R^n)}
	+
	\tfrac{1}{h_i}\int_0^t\mathrm{e}^\frac{s-t}{h_i}\|\nabla g(s)\|^q_{L^\infty(\Omega,\R^n)}\,\ds\bigg]\\
	&\le
	\bigg[ \mathrm{e}^{-\frac{t}{h_i}}\|\nabla g_o\|_{L^\infty(\Omega,\R^n)}^q
	+\big(1 -\mathrm{e}^{-\frac{t}{h_i}}\big) \| \nabla g\|_{L^\infty(\Omega_T,\R^n)}^q\bigg]\\
	&\le \max\Big\{ \|\nabla g_o\|_{L^\infty(\Omega,\R^n)}^q, \| \nabla g\|_{L^\infty(\Omega_T,\R^n)}^q\Big\}.
\end{align*}
Therefore we have
\begin{align*}
	\bigg[\biint_{\Omega_T} |\nabla g_i|^q\dx\dt\bigg]^\frac1q\le
	\max\Big\{ \|\nabla g_o\|_{L^\infty(\Omega,\R^n)}, \| \nabla g\|_{L^\infty(\Omega_T,\R^n)}\Big\}
\quad \forall\, q\ge 1.
\end{align*}
This, however,  implies the $L^\infty$-gradient bound
\begin{equation}\label{L-infty:g_i} 
	\|\nabla g_i\|_{L^\infty(\Omega_T,\R^n)}
	\le 
	\max\big\{ \|\nabla g_o\|_{L^\infty(\Omega,\R^n)}, \| \nabla g\|_{L^\infty(\Omega_T,\R^n)}
	\big\}<L.
\end{equation}
Moreover, from~\eqref{reg-gi} and \eqref{L-infty:g_i} we conclude 
\begin{equation}\label{dtDg}
	\partial_t g_i\in L^\infty\big( 0,T; W^{1,\infty}(\Omega)\big).
\end{equation}
The task now is to construct suitable comparison functions that match the regularized boundary values $g_i$ on the lateral boundary.
The difficulty here is that these must be constructed from a generic comparison function $v$ that coincides with $g$ on the lateral boundary of $\Omega_T$. Recall that a general comparison function $v\in g+L^\infty\big(0,T; W^{1,\infty}_0(\Omega)\big)$ satisfies $\partial_tv\in L^2(\Omega_T)$, $v(0)\in L^2(\Omega)$, and
$$
	\| \nabla v\|_{L^\infty(\Omega_T,\R^n)}\le\LL.
$$
We construct the comparison maps $v_i$ using the same procedure that we used to construct the regularized boundary values $g_i$, i.e.~we define $v_i:=[v]_{h_i}$. Note that the time mollification
in  \eqref{eq:time-mollification} is performed with $g_o$ instead of $v(0)$ as the initial value, i.e.
\begin{equation*}
	v_i(t):= [g]_{h_i}(t)= \mathrm{e}^{-\frac{t}{h_i}}g_{o}+\tfrac{1}{h_i}\int_0^t\mathrm{e}^\frac{s-t}{h_i}v(s)\,\ds\quad
	\mbox{for any $i\in\N$.}
\end{equation*}
The question naturally arises why $g_o$ and not $v(0)$ is used as the initial value. This is because the gradient constraint must be satisfied 
for the time mollification  $v_i$. 
This, however,  is not possible with the initial values $v(0)$, since $v(0)\in L^2(\Omega)$ but not necessarily in $W^{1,\infty}(\Omega)$. Lemma \ref{lem:moll-props} ensures that $\partial_t v_i\in L^2(\Omega_T)$. Since $v\in L^q\big( 0,T; W^{1,q}(\Omega)\big)$ and $g_o\in W^{1,q}(\Omega)$ for any $q\ge 1$, Lemma \ref{lem:moll-props} with $r=q$ and $X=W^{1,q}(\Omega)$ yields
$v_i\in L^{q}(0,T;W^{1,q}(\Omega))$, and as for the $g_i$ we obtain
 \begin{align*}
	\bigg[\biint_{\Omega_T} |\nabla v_i|^q\dx\dt\bigg]^\frac1q\le
	\max\big\{ \|\nabla g_o\|_{L^\infty(\Omega,\R^n)}, \| \nabla v\|_{L^\infty(\Omega_T,\R^n)}\big\}
\quad \forall\, q\ge 1,
\end{align*}
so that
$$
	\|\nabla v_i\|_{L^\infty(\Omega_T,\R^n)}
	\le 
	\max\big\{ \|\nabla g_o\|_{L^\infty(\Omega,\R^n)}, \| \nabla v\|_{L^\infty(\Omega_T,\R^n)}
	\big\}\le \LL.
$$
Hence, $v_i\in g_i+L^\infty \big( 0,T; W^{1,\infty}_0(\Omega)\big)$. Furthermore, by Lemma \ref{lem:moll-props}
we have
\begin{equation}\label{conv-vi}
	\mbox{$v_i\to v $ in $L^q(\Omega_T)$ in the limit $i\to\infty$}
\end{equation} 
for any $q\ge 1$.

Next we prove that the sequence $(\partial_t v_i)_{i\in\N}$ is uniformly bounded in $L^1\big( 0,T;L^2(\Omega)\big)$. For this purpose we define
\begin{equation}\label{def:tilde-vi}
	\widetilde v_i(t)
	:= 
	v_i(t)+\mathrm e^{-\frac{t}{h_i}}(v(0)-g_o)
	=
	\mathrm{e}^{-\frac{t}{h_i}}v(0)+\tfrac{1}{h_i}
	\int_0^t\mathrm{e}^\frac{s-t}{h_i}v(s)\,\ds,
\end{equation} 
so that $\widetilde v_i(0)=v(0)$. Correcting the initial values of $v_i$ from $g_o$
to $v(0)$ allows the application of the part of Lemma \ref{lem:moll-props} (with $r=2$, $X=L^2(\Omega)$)
that refers to the time derivative of the mollification,  because  $\partial_tv\in
L^2(\Omega_T)=L^2\big( 0,T, L^2(\Omega)\big)$ by assumption.
By Lemma \ref{lem:moll-props} we obtain
\begin{align*}
	\int_0^T\big\| \partial_t\widetilde v_i(t)\big\|_{L^2(\Omega)}\dt
	&\le 
	\sqrt{T}\big\| \partial_t\widetilde v_i\big\|_{L^2(\Omega_T)}
	\le
	\sqrt{T}\| \partial_tv\|_{L^2(\Omega_T)}.
\end{align*}
Moreover, we have
\begin{align*}
	\int_0^T \big\| \partial_t\big( \mathrm e^{-\frac{t}{h_i}}(v(0)-g_o)\big)\big\|_{L^2(\Omega)}\dt
	&=\int_0^T\tfrac{1}{h_i}\mathrm e^{-\frac{t}{h_i}}\dt\, \| v(0)-g_o\|_{L^2(\Omega)}\\
	&\le \| v(0)-g_o\|_{L^2(\Omega)}.
\end{align*}
This is exactly the reason why the $L^1\big(0,T; L^2(\Omega)\big)$-norm comes into play, because for any other exponents $r>1$ we only have
\begin{align*}
	\int_0^T \big\| \partial_t\big( \mathrm e^{-\frac{t}{h_i}}(v(0)-g_o)\big)\big\|_{L^2(\Omega)}^r\dt
	&=\int_0^T\tfrac{1}{h_i^r}\mathrm e^{-\frac{rt}{h_i}}\dt\, \| v(0)-g_o\|_{L^2(\Omega)}^r\\
	&= \frac{1-\mathrm e^{-\frac{rT}{h_i}}}
	{rh_i^{r-1}}\, \| v(0)-g_o\|_{L^2(\Omega)}^r\to\infty
\end{align*}
in the limit $i\to\infty$, 
which means that there is no uniform $L^r\big(0,T; L^2(\Omega)\big)$-bound with $r>1$ for $\mathrm e^{-\frac{t}{h_i}}(v(0)-g_o)$.
Together, the second and third last inequalities and \eqref{def:tilde-vi} imply
\begin{align}\label{L1-L2}
	\int_0^T\| \partial_t v_i(t)\|_{L^2(\Omega)}\dt
	&\le 
	\sqrt{T}\| \partial_tv\|_{L^2(\Omega_T)}+\| v(0)-g_o\|_{L^2(\Omega)}.
\end{align}
Finally, Lemma~\ref{lem:moll-props} applied with $r=q$ and $X=L^q(\Omega,\R^n)$ yields that 
\begin{equation}\label{Dvi-to-Dv}
	\mbox{$\nabla v_i \to \nabla v$ in $L^q(\Omega_T,\R^n)$ as $i\to\infty$}
\end{equation}
for any $q\ge 1$.
We note that \eqref{L1-L2} and \eqref{Dvi-to-Dv} apply in particular for the choice $v=g$, since $g$ is admissible as comparison map.
%

\subsection{The regularized problem formulation} According to \eqref{L-infty:g_i} and~\eqref{dtDg}, the regularized boundary values $g_i$ defined in  \eqref{def:gi} satisfy all the conditions assumed in Chapter~\ref{sec:ex-regularized}.
Hence the application of Proposition~\ref{prop:existence-reg} ensures the existence of variational solutions $u_i\in \Lip_{g_i}(\Omega_T,\LL)$ to the variational inequalities
\begin{align}\label{var-in-ui}
	\iint_{\Omega_\tau}f(\nabla u_i)\,\dx\dt
	&\le
	\iint_{\Omega_\tau}\partial_tv(v-u_i)\,\dx\dt +
    \iint_{\Omega_\tau}f(\nabla v)\,\dx\dt \nonumber\\
	&\phantom{\le\,}
	-\tfrac12\|(v-u_i)(\tau)\|^2_{L^2(\Omega)}+\tfrac12\|v(0)-g_o\|^2_{L^2(\Omega)}.
\end{align}
These  hold true for any $\tau\in[0,T]$ and any  $v\in \Lip_{g_i}(\Omega_T,\LL)$ with $\partial_tv\in L^2(\Omega_T)$ and $v(0)\in L^2(\Omega)$. 
We let $t_o\in(0,T]$ and consider $\tau\in[0,t_o]$ in \eqref{var-in-ui}. 
Choosing the testing function $v=g_i$ (recall that $g_i$ is admissible) and discarding the positive term on the right-hand side, we obtain for any $\tau\in[0,t_o]$ that
\begin{align*}
	&\tfrac12\|(g_i-u_i)(\tau)\|^2_{L^2(\Omega)}\\
	&\quad\le
    \iint_{\Omega_\tau}|\partial_t g_i| |g_i-u_i|\,\dx\dt +
    |\Omega_\tau|\sup_{|\xi|\le \LL} f(\xi) \\ 
    &\quad\le 
    \tfrac14 \sup_{t\in [0,t_o]} \|(g_i-u_i)(t)\|^2_{L^2(\Omega)} +
	t_o\|\partial_tg_i\|^2_{L^2(\Omega_{t_o})} +
    |\Omega_{t_o}|\sup_{|\xi|\le \LL} f(\xi).
\end{align*}
Taking the supremum with respect to $\tau\in[0,t_o]$, we may re-absorb the first term on the right into the left. Using the bound  \eqref{dtgi} for $\partial_tg_i$  to estimate the second term on the right, we obtain
for any $t_o\in(0,T]$ that
\begin{align}\label{gi-ui}
	\sup_{t\in [0,t_o]} \|(g_i-u_i)(t)\|^2_{L^2(\Omega)}
    \le 
	4t_o\|\partial_tg\|^2_{L^2(\Omega_{t_o})} +
    4|\Omega_{t_o}|\sup_{|\xi|\le \LL} f(\xi).
\end{align}
In particular for the choice $t_o=T$, we obtain an $L^\infty{-}L^2$ estimate for the solutions $u_i$. Indeed, we have
\begin{align*}
	\sup_{t\in [0,T]}\|u_i(t)\|^2_{L^2(\Omega)}
	\le
	2 \sup_{t\in [0,T]}\|g_i(t)\|^2_{L^2(\Omega)} +
    8T \|\partial_tg\|^2_{L^2(\Omega_T)} + 
    8|\Omega_T| \sup_{|\xi|\le \LL} f(\xi).
\end{align*}
The first term on the right-hand side is uniformly (with respect to $i\in\N$) bounded  with respect to $i\in\N$ according to Lemma \ref{lem:moll-props} applied with $r=\infty$ and
$X=L^2(\Omega)$. Therefore, the above estimate together with the uniform
gradient bound $\|\nabla u_i\|_{L^\infty (\Omega_T,\R^n)}\le\LL$ imply that the sequence $(u_i)_{i\in\N}$ is uniformly 
bounded in  $L^\infty\big( 0,T;L^2(\Omega)\big)$ and $L^\infty\big( 0,T; W^{1,\infty}(\Omega)\big)$.

\subsection{Limit passage and convergence} 
Due to the uniform estimates there exists a function 
$$
 	u\in L^{\infty}\big(0,T; L^2(\Omega)\big)\cap \bigcap_{q\ge 1}L^{q}\big(0,T;W^{1,p}(\Omega )\big)
$$
and a subsequence of $u_{i}$ (still denoted this way) such that
\begin{equation}\label{conv-ui}
	\left\{ \begin{array}{cl}
	u_{i}\rightharpoonup u & \mbox{weakly in $L^{q}\big(0,T; W^{1,q}(\Omega)\big)$ for any $q\ge 1$,}\\[5pt]
	u_{i}\rightharpoonup u & \mbox{weakly$^\ast$ in $L^{\infty}\big(0,T, L^2(\Omega)\big)$,}
\end{array}\right.
\end{equation}
as $i\to\infty$. At this point we need to show that $u$ is the sought variational solution with lateral boundary values $g$ and initial value $g_o$, that  satisfies the gradient constraint. The latter results from the weak $L^p$-convergence of the gradients and  the lower-semicontinuity of the $L^q$-norm with respect to weak convergence. More precisely, for any $q\ge 1$ we have 
$$
	\bigg[\biint_{\Omega_T} |\nabla u|^q\,\dx\dt\bigg]^\frac1q
	\le
	\liminf_{i\to\infty}
	\bigg[\biint_{\Omega_T} |\nabla u_i|^q\,\dx\dt\bigg]^\frac1q
	\le \LL,
$$
so that the gradient constraint $\|\nabla u\|_{L^\infty (\Omega_T,\R^n)}\le\LL$ holds. In order to check that $u$ solves the variational inequality we consider an arbitrary comparison function $v\in \Lip_g(\Omega_T,\LL)$ with $\partial_tv\in L^2(\Omega)$ and $v(0)\in L^2(\Omega)$. As before, we write $v_i:=[v]_{h_i}$ for the time mollification of $v$  and with initial value $g_o$. As shown in \S\,\ref{sec:constr-data}, the $v_i$ are admissible in the variational inequality~\eqref{var-in-ui}, from which we obtain
\begin{align}\label{var-ungl}
	\tfrac12 \| (v_i-u_i)(\tau)\|^2_{L^2(\Omega)} &+\iint_{\Omega_\tau} f(\nabla u_i)\,\dx\dt\nonumber\\
	&\le \iint_{\Omega_\tau}\partial_tv_i(v_i-u_i)\,\dx\dt +
    \iint_{\Omega_\tau} f(\nabla v_i)\,\dx\dt
\end{align}
for any $\tau\in [0,T]$. Note that $v_i$ and $g_i$ attend the same initial datum $g_o$, so that the $L^2$-boundary term vanishes for $t=0$.
Now, the overall goal is to pass in \eqref{var-ungl} to the limit $i\to\infty$. Hereby, the degree of difficulty in the treatment of the individual terms is quite different. We start with the second integral on the right. Since both functions $v$ and $v_i$ satisfy the gradient constraint
we have
\begin{align}\label{cov-f-vi}
	\bigg|\iint_{\Omega_\tau} f(\nabla v_i)\,\dx\dt &
    -\iint_{\Omega_\tau} f(\nabla v)\,\dx\dt\bigg|\nonumber\\
	&\le 
 \sup_{|\xi|\le\LL} |\nabla f(\xi)|\iint_{\Omega_T} | \nabla v_i-\nabla v|\,\dx\dt.
\end{align}
Because of \eqref{Dvi-to-Dv} with $q=1$, the integral on the right converges to 0. This allows us to replace the second integral on the right-hand
side of \eqref{var-ungl}  by $\iint_{\Omega_\tau} f(\nabla v)\,\dx\dt$ after passing to the limit $i\to\infty$. In the first integral, i.e.~the integral containing the time derivative, we replace $\partial_t v_i$ by $\partial_t\widetilde v_i$, where $\widetilde v_i$ is defined in \eqref{def:tilde-vi}. 
We recall that passing from $v_i$  to $\widetilde v_i$ means just a correction of the initial values. We have
\begin{align}\label{time-int}
	\iint_{\Omega_\tau}&\partial_tv_i(v_i-u_i)\,\dx\dt\nonumber\\
	&=
	\iint_{\Omega_\tau}\partial_t\widetilde v_i(v_i-u_i)\,\dx\dt
	+
	\tfrac1{h_i} \iint_{\Omega_\tau}\mathrm e^{-\frac{t}{h_i}}(v(0)-g_o)(v_i-u_i)\,\dx\dt.
\end{align}
In view of the second part of Lemma \ref{lem:moll-props} we have $\partial_t\widetilde v_i(t) = \frac{1}{h_i}\int_0^t\mathrm e^\frac{s-t}{h_i}\partial_sv(s)\,\ds$. 
Thus the first part of Lemma \ref{lem:moll-props} can be applied to $\partial_tv\in L^2(\Omega_T)=L^2\big( 0,T;L^2(\Omega)\big)$ with the choice $v_o=0$. This yields $\partial_t\widetilde v_i\to \partial_tv$ in $L^2(\Omega_T)$ as $i\to\infty$. Combining this observation with \eqref{conv-vi} and
\eqref{conv-ui}$_2$, we get for the first term on the left
\begin{equation}
	\lim_{i\to\infty}
	\iint_{\Omega_\tau}\partial_t\widetilde v_i(v_i-u_i)\,\dx\dt
	=
	\iint_{\Omega_\tau}\partial_t v(v-u)\,\dx\dt.
\end{equation}
Next, we write the second integral on the right-hand side of \eqref{time-int} in the form
\begin{align*}
	\tfrac1{h_i} \iint_{\Omega_\tau}\mathrm e^{-\frac{t}{h_i}}(v(0)-g_o)(v_i-u_i)\,\dx\dt
	&=
	\mathbf{I}_i+\mathbf{II}_i+\mathbf{III}_i+\mathbf{IV}_i,
\end{align*}
where we abbreviated
\begin{align*}
	\mathbf{I}_i
	&
	:=\tfrac1{h_i}\iint_{\Omega_\tau}\Big(\mathrm e^{-\frac{2t}{h_i}}- \mathrm e^{-\frac{t}{h_i}}\Big) (v(0)-g_o)g_o\,\dx\dt,\\
	\mathbf{II}_i
	&
	:=\tfrac1{h_i^2}\iint_{\Omega_\tau}\bigg[\mathrm e^{-\frac{2t}{h_i}} (v(0)-g_o)\int_0^t \mathrm e^\frac{s}{h_i}v(0)\ds\bigg]\,\dx\dt,\\
	\mathbf{III}_i
	&
	:=\tfrac1{h_i^2}\iint_{\Omega_\tau}\bigg[\mathrm e^{-\frac{2t}{h_i}} (v(0)-g_o)\int_0^t \mathrm e^\frac{s}{h_i}(v(s)-v(0))\ds\bigg]\,\dx\dt,\\
	\mathbf{IV}_i
	&
	:=\tfrac1{h_i}\iint_{\Omega_\tau} \mathrm e^{-\frac{t}{h_i}} (v(0)-g_o)(g_o-u_i)\,\dx\dt.
\end{align*}
In $\mathbf{I}_i$ we first compute the one-dimesional integral with respect to $t$ and then pass to the limit  $i\to\infty$.
This leads to
\begin{align*}
	\lim_{i\to\infty}\mathbf{I}_i
	&=
	\lim_{i\to\infty} \tfrac1{h_i}\int_0^\tau \Big(\mathrm e^{-\frac{2t}{h_i}}- \mathrm e^{-\frac{t}{h_i}}\Big) \,\dt
	\int_\Omega (v(0)-g_o)g_o\,\dx\\
	&=-\tfrac12 \int_\Omega (v(0)-g_o)g_o\,\dx.
\end{align*}
The second integral $\mathbf{II}_i$ is treated analogously. We obtain
\begin{align*}
	\lim_{i\to\infty}\mathbf{II}_i
	&=
	\lim_{i\to\infty} \tfrac1{h_i^2}\int_0^\tau  \mathrm e^{-\frac{2t}{h_i}}
	\int_0^t\mathrm e^{\frac{s}{h_i}}\,\ds\dt
	\int_\Omega (v(0)-g_o)v(0)\,\dx\\
	&=\tfrac12 \int_\Omega (v(0)-g_o)v(0)\,\dx.
\end{align*}
Together with the result for $\mathbf{I}_i$ this gives
\begin{align*}
	\lim_{i\to\infty}\big[ \mathbf{I}_i+\mathbf{II}_i\big]
	&=\tfrac12 \| v(0)-g_o\|^2_{L^2(\Omega)},
\end{align*}
which is exactly the $L^2(\Omega)$-boundary term at $t=0$. Next, we treat the term $\mathbf{III}_i$.
This integral can can be estimated with H\"older's inequality. We obtain
\begin{align*}
	\lim_{i\to\infty}&\mathbf{III}_i\\
	&\le
	\| v(0)-g_o\|_{L^2(\Omega)}
	\lim_{i\to\infty}\tfrac1{h_i^2}\int_0^\tau
	\mathrm e^{-\frac{2t}{h_i}} \int_0^t \mathrm e^\frac{s}{h_i}
	\bigg[ \int_\Omega |v(s)-v(0)|^2\,\dx \bigg]^\frac12\,\ds\dt\\
	&=
	\| v(0)-g_o\|_{L^2(\Omega)}
	\lim_{i\to\infty}\tfrac1{h_i^2}\int_0^\tau
	\mathrm e^{-\frac{2t}{h_i}} \int_0^t \mathrm e^\frac{s}{h_i}
	\bigg[ \int_\Omega \bigg|\int_0^s\partial_\sigma v(\sigma)\dsigma\bigg|^2\,\dx \bigg]^\frac12\,\ds\dt\\
	&\le
	\| v(0)-g_o\|_{L^2(\Omega)}\| \partial_tv\|_{L^2(\Omega_\tau)}
	\lim_{i\to\infty}\tfrac1{h_i^2}\int_0^\tau \mathrm e^{-\frac{2t}{h_i}}\int_0^t\sqrt{s}\mathrm e^\frac{s}
	{h_i}\,\ds\dt\\
	&\le
	\tfrac23
	\| v(0)-g_o\|_{L^2(\Omega)}\| \partial_tv\|_{L^2(\Omega_\tau)}
	\lim_{i\to\infty}\tfrac1{h_i^2}
	\int_0^\tau \mathrm e^{-\frac{t}{h_i}}t^\frac32\,\dt\\
	&\le
	\tfrac23
	\| v(0)-g_o\|_{L^2(\Omega)}\| \partial_tv\|_{L^2(\Omega_\tau)}
	\int_0^\infty s^\frac32\mathrm e^{-s}\,\ds
	\lim_{i\to\infty}\sqrt{h_i} 
    =
    0.
\end{align*}
Finally, we deal with term $\mathbf{IV}_i$, which is the most difficult one. First, we apply
the Cauchy-Schwarz inequality in space. In the resulting integral we decompose the time interval $[0,\tau]$
into two parts, i.e.~in $[0,\sqrt{h_i}]$ and $(\sqrt{h_i},\tau]$. This leads to
\begin{align*}
	| \mathbf{IV}_i|
	&\le 
	\| v(0)-g_o\|_{L^2(\Omega)}
	\tfrac{1}{h_i}\int_0^\tau \mathrm e^{-\frac{t}{h_i}}\|g_o-u_i(t)\|_{L^2(\Omega)}\,\dt\\
	&=
    \| v(0)-g_o\|_{L^2(\Omega)}
	\big[ 
	 \mathbf{IV}_i^{(1)} +\mathbf{IV}_i^{(2)}
	\big],
\end{align*}
where we used the shorthand notation 
\begin{align*}
	\mathbf{IV}_i^{(1)}
	&:=\tfrac{1}{h_i}\int_0^{\sqrt{h_i}} \mathrm e^{-\frac{t}{h_i}}\|g_o-u_i(t)\|_{L^2(\Omega)}\,\dt,\\
	\mathbf{IV}_i^{(2)}
	&:=
	\tfrac{1}{h_i}\int_{\sqrt{h_i}}^\tau \mathrm e^{-\frac{t}{h_i}}\|g_o-u_i(t)\|_{L^2(\Omega)}\,\dt.
\end{align*}
Here the second term $\mathbf{IV}_i^{(2)}$ is easier to deal with, and therefore we start with it.
Using the uniform $L^\infty (0,T;L^2(\Omega))$-bound for $u_i$ we obtain
\begin{align*}
	\lim_{i\to\infty}\mathbf{IV}_i^{(2)}
	&\le
	\lim_{i\to\infty}\|g_o-u_i\|_{L^\infty (0,T;L^2(\Omega))}
	\tfrac{1}{h_i}\int_{\sqrt{h_i}}^\tau \mathrm e^{-\frac{t}{h_i}}\,\dt\\
	&\le
	\sup_{i\in\N} \|g_o-u_i\|_{L^\infty (0,T;L^2(\Omega))}
	\lim_{i\to\infty}
	\Big[ \mathrm e^{-\frac{1}{\sqrt{h_i}}}-\mathrm e^{-\frac{\tau}{h_i}}\Big]\\
	&\le 
	\sup_{i\in\N} \|g_o-u_i\|_{L^\infty (0,T;L^2(\Omega))}
	\lim_{i\to\infty}
	 \mathrm e^{-\frac{1}{\sqrt{h_i}}}
    =0.
\end{align*}
Thus, it remains to consider the term $\mathbf{IV}_i^{(1)}$. We have
\begin{align*}
	\mathbf{IV}_i^{(1)}
	&\le
	\| g_o-u_i\|_{L^\infty(0,\sqrt{h_i};L^2(\Omega))}
	\tfrac{1}{h_i}\int_0^{\sqrt{h_i}} \mathrm e^{-\frac{t}{h_i}}\,\dt\\
	&=
	\| g_o-u_i\|_{L^\infty(0,\sqrt{h_i};L^2(\Omega))}
	\Big( 1-\mathrm e^{-\frac{1}{\sqrt{h_i}}}\Big)\\
	 &\le
	\| g_o-u_i\|_{L^\infty(0,\sqrt{h_i};L^2(\Omega))}.
\end{align*}
From \eqref{gi-ui} applied with $t_o=\sqrt{h_i}$ we obtain
\begin{align*}
	\lim_{i\to\infty} \sup_{t\in [0,\sqrt{h_i}]} \|(g_i-u_i)(t)\|^2_{L^2(\Omega)}
    \le 
	4\lim_{i\to\infty} \sqrt{h_i} \Big[\|\partial_tg\|^2_{L^2(\Omega_{T})} +
    |\Omega|\sup_{|\xi|\le \LL} f(\xi)\Big]
    =0.
\end{align*}
At this point, it remains to analyze the difference between $g_i(t)$ and $g_o$ close to $t=0$.
By~\eqref{def:gi} we have for any $t\in [0,\sqrt{h_i}]$ that 
\begin{align*}
	\| g_i(t)-g_o\|^2_{L^2(\Omega)}
	&= 
	\int_{\Omega}\bigg|\big(\mathrm e^{-\frac{t}{h_i}}-1\big) g_o + 
    \tfrac{1}{h_i}\int_0^t\mathrm{e}^\frac{s-t}{h_i}g(s)\,\ds \bigg|^2 \,\dx \\
	&= 
	\int_{\Omega}\bigg|
    \tfrac{1}{h_i}\int_0^t\mathrm{e}^\frac{s-t}{h_i}(g(s)-g_o)\,\ds \bigg|^2 \,\dx \\
    &\le 
    \sup_{\tau\in [0,\sqrt{h_i}]} \| g(\tau)-g_o \|^2_{L^2(\Omega)}.
\end{align*}
Taking the supremum with respect to $t$ over the interval $[0,\sqrt{h_i}]$, we find
\begin{align*}
	\sup_{t\in [0,\sqrt{h_i}]} \| g_i(t)-g_o\|^2_{L^2(\Omega)}
	&\le 
	\sup_{t\in [0,\sqrt{h_i}]} \| g(t)-g_o \|^2_{L^2(\Omega)}.
\end{align*}
Since $g\in C^0\big( [0,T]; L^2(\Omega)\big)$ with $g(0)=g_o$, the right-hand side of the preceding
inequality converges to $0$ as $i\to\infty$. The combination of the last inequalities implies that also $\mathbf{IV}_i^{(1)}\to 0$ in the limit $i\to\infty$. Therefore we have everything at hand to pass to the limit in \eqref{time-int}. 
For a.e.~$\tau\in (0,T)$ we get
\begin{align}\label{time-int-limit}
	\lim_{i\to\infty}
	\iint_{\Omega_\tau}&\partial_tv_i(v_i-u_i)\,\dx\dt
	=
	\iint_{\Omega_\tau}\partial_t v(v-u)\,\dx\dt
	+
	\tfrac12 \| v(0)-g_o\|^2_{L^2(\Omega)}.
\end{align}
Moreover, the first integral on the right-hand side can be uniformly bounded with respect to $\tau\in (0,T)$ and $i\in\N$ using \eqref{L1-L2}.
Indeed, we have
\begin{align}\label{time-int-bd}
	\sup_{i\in\N}&\sup_{\tau\in (0,T)}\bigg|\iint_{\Omega_\tau}\partial_tv_i(v_i-u_i)\,\dx\dt\bigg|\nonumber\\
	&\le
	\sup_{i\in\N}\|\partial_tv_i\|_{L^1(0,T;L^2(\Omega))}\|v_i-u_i\|_{L^\infty(0,T;L^2(\Omega))}\nonumber\\
	&\le
	\sqrt{T}\|\partial_tv\|_{L^2(\Omega_T)}\sup_{i\in\N}\|v_i-u_i\|_{L^\infty(0,T;L^2(\Omega))}<\infty.
\end{align}
Before proceeding to the limit $i\to\infty$ in \eqref{var-ungl}, we integrate both sides with respect to $\tau$ over time intervals $(t_o, t_o+\delta)$ contained in $(0,T)$ and then take means on both sides. In this way we get
\begin{align}\label{var-ungl-mittel}
	\tfrac12 \mint_{t_o}^{t_o+\delta}&\| (v_i-u_i)(\tau)\|^2_{L^2(\Omega)}\,\dtau +
	\iint_{\Omega_{t_o}} f(\nabla u_i)\,\dx\dt \nonumber\\
	&\le 
    \mint_{t_o}^{t_o+\delta}\bigg[ \iint_{\Omega_\tau} \partial_tv_i(v_i-u_i)\,\dx\dt +\iint_{\Omega_\tau} f(\nabla v_i)\,\dx\dt\bigg]\dtau 
\end{align}
The integrals on the left-hand side are both lower semicontinuous. Note that according to \eqref{conv-vi} $v_i$ converges strongly in $L^2(\Omega)$ to $v$, that $f$ is convex, and that according to \eqref{conv-ui}$_1$ $\nabla u_i$ converges weakly in $L^q(\Omega_T,\R^n)$ to $\nabla u$ for any $q\ge 1$. According to \eqref{cov-f-vi} the integral of $f(\nabla v_i)$ converges to the corresponding integral of $f(\nabla v)$. Finally, by the dominated convergence theorem we can also pass to the limit $i\to\infty$ in the integral on the right-hand side containing the time derivative. Recall that for a.e.~$\tau$ the functions   $\tau\mapsto \iint_{\Omega_\tau}\partial_tv_i(v_i-u_i)\,\dx\dt$ are converging, i.e.~on the one hand we have a.e.~convergence on $(0,T)$, and on the other hand,  the functions are uniformly bounded with respect to $i$ and $\tau$. In this way we get
\begin{align*}
	\tfrac12 \mint_{t_o}^{t_o+\delta}&
    \| (v-u)(\tau)\|^2_{L^2(\Omega)}\,\dtau +
	\iint_{\Omega_{t_o}} f(\nabla u)\,\dx\dt \nonumber\\
	&\le \mint_{t_o}^{t_o+\delta}\bigg[ \iint_{\Omega_\tau}\big[\partial_tv(v-u)+f(\nabla v)\big]\,\dx\dt\bigg]\dtau 
	+\tfrac12 \| v(0)-g_o\|^2_{L^2(\Omega)}
\end{align*}
for any $t_o\in[0,T]$ and any $\delta\in (0,T-t_o]$. Here we let $\delta\downarrow 0$ and obtain
\begin{align*}
	\tfrac12 \| (v&-u)(t_o)\|^2_{L^2(\Omega)} +
	\iint_{\Omega_{t_o}} f(\nabla u)\,\dx\dt \nonumber\\
	&\le \iint_{\Omega_{t_o}}\big[\partial_tv(v-u)+f(\nabla v)\big]\dx\dt
	+\tfrac12 \| v(0)-g_o\|^2_{L^2(\Omega)}
\end{align*}
for a.e.~$t_o\in[0,T]$. Overall, we have  shown that $u$ is a variational solution of the gradient constraint problem. The uniqueness can be deduced as in~\cite[Lemma~3.3]{BDS}. This finishes the proof of Proposition~\ref{prop:ex-constr}.

\section{Eliminating the gradient constraint}\label{sec:remove}

In this section we will remove the gradient constraint from the variational solutions constructed in Proposition~\ref{prop:ex-constr}. This will be achieved by showing that the gradient is bounded in terms of a constant depending only on the data. Throughout this section we assume all hypothesis of Theorem~\ref{thm:existence-var-sol} to be in force.

\subsection{Properties of the integrand function $f$}\label{sec:assump-bds}

The following proposition contains some useful properties of  $f$ and its polar function (convex conjugate) $f^\ast\colon \R^n\to [-\infty,\infty]$ defined by
\begin{equation*}
    f^\ast(\eta):=\sup_{\xi \in\R^n}\big\{ \eta \cdot\xi-f(\xi)\big\}.
\end{equation*}

\begin{prop}\label{propertyf}
Suppose that $f\colon \R^n\rightarrow\R$ satisfies the set of assumptions~\ref{A2}. Then, we have
\begin{enumerate}[label=\textnormal{(\arabic*)}, ref=(\arabic*), wide=0pt, leftmargin=*, topsep=7pt]
\item there exists $c=c(\epsilon)\in\R$ such that
\begin{equation*}
    f(\xi)\ge \tfrac14 \epsilon|\xi|^2 - c  \quad\mbox{for any $\xi\in\R^n$;}
\end{equation*}
\item $f^*$ is superlinear and at most of quadratic qrowth, i.e.~that
\begin{equation*}
    f^*(\eta)\le \tfrac{2}{\epsilon}|\eta|^2 + c  \quad\mbox{for any $ \eta\in\R^n$;}
\end{equation*}
\item there exists $r>0$ such that $f^*\in C^{1,1}(\R^n\setminus B_r(0))$ and 
$$
    |D^2 f^*(\eta)|\le\tfrac{1}{\epsilon}
    \quad\mbox{for a.e.~$\eta\in \R^n\setminus B_r(0)$;}
$$
\item the restriction of $\nabla f$ to $\R^n\setminus B_1$  is an invertible vector-field with $$(\nabla f)^{-1}=\nabla f^\ast.$$ 
\end{enumerate}
 \end{prop}
\begin{proof}[Proof]
Property (1) follows for large values of $|\xi|$ by Taylor's expansion and the uniform convexity of $f$ from~\eqref{hess}. Property (2) is a straightforward consequence of the definition of $f^*$ taking into account that $f$ is super-linear (cf.~\cite[Lemma 3.1]{Giannetti-Treu:2022})
and the assertion from (1). The proof of property (3) can be found in \cite{Fiaschi-Treu:2012} for a special class of functions. The general case can be retrieved  from \cite[Lemma 3.3]{Giannetti-Treu:2022}. 
Property (4) is a consequence of \cite[Chapter 1, Corollary 5.2]{Ekeland} and  the fact that $f$ is a $C^2$-function outside the unit ball by assumption.
\end{proof}
 

\subsection{Construction of barriers}

Let $\alpha\in\R\setminus\{0\}$. For $y\colon[0,T]\to\R^n$ and $c\colon[0,T]\to\R$ we consider the following function
\begin{equation}\label{def:barrier}
    v(x,t)
    =
    \frac{n}{\alpha}f^*\Big(\frac{\alpha}{n}\big(x-y(t)\big)\Big)-c(t).
\end{equation}
In view of Proposition~\ref{propertyf}\,(3),\,(4) it is easy to check that 
\begin{equation}\label{cellina}
    \Div_x \nabla_\xi f\big(\nabla v(x,t)\big)=\alpha,
\end{equation}
provided that 
\begin{equation*}
    \Big|\frac{\alpha}{n}\big(x-y(t)\big)\Big| \ge r
    \quad\mbox{and}\quad 
    \Big|\nabla f^*\Big(\frac{\alpha}{n}\big(x-y(t)\big)\Big)\Big|
    \ge 1.
\end{equation*}
We recall that the class of functions in \eqref{def:barrier} was introduced by Cellina in \cite{Cellina:2007} for the stationary case. In the same paper it was also shown that this class of functions has a property corresponding to \eqref{cellina} at the level of the associated variational integral, i.e.~they minimize the variational integral that formally has \eqref{cellina} as Euler-Lagrange equation. For this, it was only necessary to assume that $f$ is convex. 
Our goal is to show that for  suitable
choices of $\alpha$, $y(t)$ and $c(t)$ the function $v$ can be used in the time dependent setting as a barrier function. We start with an auxiliary lemma.

\begin{lem}\label{lem:b}
Let $\alpha\in\R\setminus\{0\}$ and $M,Q_1>0$. 
For $b\in \R$ and $w\in \R^n$ we define
$$
	\widetilde\Omega_{b,w}
	:=
	\left\{ x\in\R^n: \frac{n}{\alpha}f^*\left(\frac{\alpha}{n}\,x 
 \right)-w\cdot x-b\le 0\right\}.
$$
If 
$$
	b
	>
	\min_{x\in\R^n} 
	\Big[ \frac{n}{\alpha}f^*\left(\frac{\alpha}{n}\,x 
 \right)-w\cdot x\Big],  
$$
then $\widetilde\Omega_{b,w}$ is non-empty, bounded and convex. Furthermore, if
\begin{equation}\label{choice-b}
	b
	\ge
	\Gamma
	:=
	\max_{x\in\overline{B_M(0)},\,w\in\overline{B_{Q_1}(0)} }\Big[ \frac{n}{\alpha}f^*\left(\frac{\alpha}{n}\,x 
 \right)-w\cdot x\Big],
\end{equation}
then, for any $w\in\overline{B_{Q_1}(0)}$ the ball $\overline{B_M(0)}$ is contained in $\widetilde\Omega_{b,w}$. 
\hfill $\Box$ 

\end{lem}

\begin{proof}[Proof]
The proof is straightforward. For the sake of completeness, however, we will elaborate on it.
Since $f^\ast$ has superlinear growth and $\R^n\ni x\mapsto w\cdot x$ is linear, we have 
$$
    k(x):=\frac{n}{\alpha}f^*\left(\frac{\alpha}{n}\,x 
    \right)-w\cdot x\to\infty\quad\mbox{as $|x|\to\infty$.}
$$
This implies that $k$ attains its minimum on $\R^n$. If $b$ is larger than the minimal value of $k$, we obtain that 
$\widetilde \Omega_{b,w}$ is non-empty and bounded. In addition, $\widetilde \Omega_{b,w}$ is convex because it is constructed as a sublevel set of a convex function. This proves the first assertion of the lemma.

The second assertion, i.e.~that \eqref{choice-b} implies $\overline{B_M(0)}\subset \widetilde \Omega_{b,w}$, follows similarly. Fix $w\in \overline{B_{Q_1}(0)}$. Then $k$ is bounded on $\overline{B_M(0)}$. If $b$ is greater than the maximum value of $k$ on $\overline{B_M(0)}$, then $\overline{B_M(0)}$ is a subset of  $\widetilde \Omega_{b,w}$, proving the second claim.
\end{proof}

\begin{prop}\label{prop:barrier}
Let the assumptions of Theorem~\ref{thm:existence-var-sol} be in force. Then, there exists $\alpha_o>0$ depending only on $n,\epsilon, \diam(\Omega), \|D^2f\|_{L^\infty(\R^n\setminus B_1)}, \|\partial_t g\|_\infty, \mathsf Q$ such that for any $\alpha\in\R$ with $|\alpha|\ge \alpha_o$ and any $x_o\in\partial\Omega$ the functions $y(t)$ and $c(t)$ can be chosen in such a way that the function $v$ defined in \eqref{def:barrier} satisfies the following properties:
\begin{enumerate}[label=\textnormal{(\arabic*)}, ref=(\arabic*), wide=0pt, leftmargin=*, topsep=7pt]
    \item $v(x_o,t)=g(x_o,t)$ for any $t\in [0,T)$;
    \item $v(x,t)\le g(x,t)$ for any $x\in\Omega$ and $t\in [0,T)$ if $\alpha>0$; if $\alpha<0$ the converse inequality holds true;
    \item if $\alpha>0$ is large enough, then $v$ is a sub-solution to the parabolic equation \eqref{eq:strong-Cauchy-Dirichlet-prob}; similarly,
    if $\alpha<0$ is chosen small enough, then
    $v$ is a super-solution.
    \item $v\in \Lip(\Omega_T)$ with  Lipschitz constant $L$ depending only on
    the data $n$, $\epsilon$, $R$, $\diam(\Omega)$, $f$, $\nabla f$, $Q$, 
    $\|\nabla g\|_{L^\infty (\Omega_T,\R^n)}$, and $\alpha$.
\end{enumerate}
\end{prop}

Note that for the purposes of what follows, it would be sufficient for condition (2) to hold on the parabolic boundary $\partial_{\mathcal{P}}\Omega_T$. However, the proof establishes the stronger result that condition (2) is satisfied throughout the entire cylinder $\Omega_T$.

\begin{proof}[Proof]
We only prove the assertions for positive $\alpha$, since the proof for the other case is analogous. 
We may assume $\alpha\ge \alpha_o\ge n$ by choosing 
$\alpha_o\ge n$. 

\textit{Step 1: An auxiliary set.}
Let $R$ be the radius from the $R$-uniformily convexity assumption~\ref{A1} on $\Omega$ and $\nu_{x_o}$ the outward pointing unit vector associated with $x_o\in\partial\Omega$; see Definition~\ref{def:unif-conv}. From Remark~\ref{rem:unif-conv} we conclude that $\Omega\subset B_R(x_o-R\nu_{x_o})$. 

Next, we let $w_{x_o}^-\colon[0,T]\to\R^n$ be the function from the $t-{\rm BSC}_Q$ condition in~\ref{A3}; see Definition~\ref{def:bdslope}. 
By $\widetilde w_{x_o}^- \colon [0,T]\to \R^n$ we denote the modified function constructed in Lemma~\ref{lem:tBSC-G} satisfying $\|\widetilde w_{x_o}^-\|_{L^\infty ([0,T],\R^n)}\le Q_1:=Q+ \|\nabla g\|_{L^\infty (\Omega_T,\R^n)}$ and $\|(\widetilde w_{x_o}^-)'\|_{L^\infty ([0,T],\R^n)} = \|(w_{x_o}^-)'\|_{L^\infty ([0,T],\R^n)}\le \mathsf Q$
and such that
\begin{align}\label{g-tilde}
    \widetilde g(x,t)
    :=
    g(x_o,t)+\widetilde w_{x_o}^{-}(t)\cdot (x-x_o)
    \le
    g(x,t)
\end{align}
for any $(x,t)\in \overline\Omega\times [0,T)$.
For any $t\in [0,T)$ we define the set 
\begin{equation*}
    \widetilde\Omega_t
    :=
    \bigg\{x\in\R^n: \underbrace{\frac{n}{\alpha}f^*\Big(\frac{\alpha}{n}\big(x-y(t)\big)\Big)-c(t)}_{=\, v(x,t)}\le \widetilde g(x,t)\bigg\}.
\end{equation*}

\textit{Step 2: Proof of claim~(1).} Here we choose $y(t)$ and $c(t)$ in such a way that 
\begin{equation}\label{step2}
    x_o\in\partial \widetilde\Omega_t
    \quad\mbox{for any $t\in [0,T)$.}
\end{equation}
In view of the definition of $\widetilde\Omega_t$ this implies claim~(1).  
%

By $\epsilon$ and $r$ we denote the corresponding parameters from~\eqref{hess} and Proposition~\ref{propertyf}\,(3).  
Due to the super-linearity of $f^\ast$ there exists $M\ge r+\diam(\Omega)$, depending only on $\epsilon, R, \diam(\Omega), \nabla f, Q_1$, such that 
\begin{equation}\label{nfs}
    |\nabla f^\ast(\eta)|>\frac{R}{\epsilon}+Q_1\qquad 
    \mbox{for any $\eta\in\R^n\setminus B_{M}(0)$.}
\end{equation}
We define $\Gamma$ according to \eqref{choice-b} from Lemma~\ref{lem:b} for this particular choice of $M$. Note that $\Gamma$ depends on $n, \epsilon, R, \diam(\Omega), \nabla f, Q_1, \alpha$.
Again, due to the super-linearity of $f^\ast$ and the upper bound $|\widetilde w_{x_o}^-(t)|\le Q_1$ for any $t\in[0,T)$ there exists $\varrho_o\ge M$ such that 
\begin{equation}\label{bgeg}
    B(t,\eta)
    :=
    \frac{n}{\alpha}f^*\left(\frac{\alpha}{n}\,\eta
 \right) - \widetilde w_{x_o}^-(t)\cdot \eta
 \ge 
 \Gamma
 \quad\mbox{for any $\eta\in\R^n\setminus B_{\varrho_o}$ and $t\in[0,T)$.}
\end{equation}
Note that $\varrho_o$ can be chosen in dependence on $n, \epsilon, R, \diam(\Omega), f, \nabla f, Q_1, \alpha$. Next we fix $\lambda>Q_1+\max\{1, r\}$ such that
\begin{equation}\label{zlambdaesplicit}
 z_\lambda(t):=\frac{n}{\alpha}\nabla f\big(\widetilde w_{x_o}^-(t)+\lambda\nu_{x_o}\big)
 \in 
 \R^n\setminus B_{\varrho_o} 
 \quad\mbox{for any $t\in[0,T)$.}
\end{equation}
This can be achieved,  because  $|\widetilde w_{x_o}^-(t)|\le Q_1$ and since $ f$ has super-linear growth. Note that $\lambda$ can be chosen in dependence on $n, \epsilon, R, \diam(\Omega), f, \nabla f, Q_1, \alpha$. Moreover, $|\widetilde w_{x_o}^-(t)+\lambda\nu_{x_o}|\ge \lambda -Q_1>1$, which ensures that the gradient in~\eqref{zlambdaesplicit} is well defined. 
Using Lemma~\ref{propertyf}\,(4), i.e.~the fact that $(\nabla f)^{-1}=\nabla f^\ast$ on $\R^n\setminus B_1$, \eqref{zlambdaesplicit} can be re-written in the form 
\begin{equation}\label{zlambdaimplicit}
 \nabla f^\ast\left(\frac{\alpha}{n}\,z_\lambda(t)\right)-\widetilde w_{x_o}^-(t)=\lambda\nu_{x_o}.
\end{equation}
Note that $|\frac{\alpha}{n}z_\lambda(t)|\ge |z_\lambda(t)|\ge \varrho_o\ge M>r$, so that by Lemma~\ref{propertyf}\,(3), $\nabla f^\ast\left(\frac{\alpha}{n}\,z_\lambda(t)\right)$ is well defined. 
Next, we let
\begin{equation*}
    b(t)
    :=
    B(t,z_\lambda(t))
    =
    \frac{n}{\alpha}f^*\left(\frac{\alpha}{n}\,z_\lambda(t) 
 \right) - \widetilde w_{x_o}^-(t)\cdot z_\lambda(t)
 \quad\mbox{for $t\in[0,T)$}
\end{equation*}
and observe that \eqref{zlambdaesplicit} and~\eqref{bgeg} imply
\begin{equation}\label{bgeg-}
    b(t)
    \ge
    \Gamma 
    \quad\mbox{for any $t\in[0,T)$.}
\end{equation}
For any $t\in[0,T)$ we now define the set
\begin{equation}\label{def:om-b}
    \widetilde\Omega_{b(t),\widetilde w_{x_o}^-(t)}=\Big\{ x\in\R^n\colon
    \frac{n}{\alpha}f^*\left(\frac{\alpha}{n}\,x 
    \right)-\widetilde w_{x_o}^-(t)\cdot x-b(t)\le 0
 \Big\}.
\end{equation}
By definition we have  $z_\lambda(t)\in\partial\widetilde\Omega_{b(t),\widetilde w_{x_o}^-(t)}$ for any $t\in[0,T)$. 
Moreover, \eqref{zlambdaimplicit} implies that the sets $\widetilde\Omega_{b(t),\widetilde w_{x_o}^-(t)}$ have outward normal $\nu_{x_o}$ at  $z_\lambda(t)\in\partial\widetilde\Omega_{b(t),\widetilde w_{x_o}^-(t)}$ for any $t\in[0,T)$. 
Finally, we define
\begin{equation}\label{def:y(t)}
    y(t):=x_o-z_\lambda(t),
\end{equation}
and 
\begin{align}\label{c(t)}
  c(t)
  &:=
  b(t) + 
  \widetilde w_{x_o}^-(t)\cdot z_\lambda(t) -g(x_o,t) \nonumber\\
  &\ =
  \frac{n}{\alpha}f^*\left(\frac{\alpha}{n}\big(x_o-y(t)\big)\right)
  -g(x_o,t)
\end{align}
and observe that
$$
\widetilde\Omega_t = y(t)+\widetilde\Omega_{b(t),\widetilde w_{x_o}^-(t)}.
$$
Since $z_\lambda(t)\in\partial\widetilde\Omega_{b(t),\widetilde w_{x_o}^-(t)}$, we have $x_o\in\partial \widetilde\Omega_t$ for any $t\in[0,T)$, which shows~\eqref{step2}. Moreover, $\widetilde\Omega_t$ has outward normal $\nu_{x_o}$ at $x_o\in\partial \widetilde\Omega_t$ for any $t\in[0,T)$. 
\medskip 

\textit{Step 3: Proof of claim~(2).}
Here we prove that 
\begin{equation}\label{step3}
    \Omega\subset\widetilde\Omega_t
    \quad\mbox{for any $t\in [0,T)$.}
\end{equation}
Due to of the definition of $\widetilde\Omega_t$ and \eqref{g-tilde} this implies claim~(2). 

We fix a time $t\in[0,T)$ and abbreviate $b=b(t)$ and $w=\widetilde w_{x_o}^-(t)$. With these abbreviations the set in \eqref{def:om-b} can be re-written as 
\begin{equation*}
    \widetilde\Omega_{b,w}=\Big\{ x\in\R^n\colon
    \frac{n}{\alpha}f^*\left(\frac{\alpha}{n}\,x 
    \right)-w\cdot x-b\le 0
 \Big\}.
\end{equation*}
Since $|w|\le Q_1$ and $b\ge\Gamma$ by~\eqref{bgeg-}, Lemma~\ref{lem:b} ensures that $\overline{B_{M}(0)}\subset \widetilde\Omega_{b,w}$, so that $|x|\ge M\ge r+\diam(\Omega)$ for any $x\in \partial \widetilde\Omega_{b,w}$. Noting that  $\alpha\ge  \alpha_o\ge n$ we conclude 
\begin{equation}\label{boundary-point}
	\big|\tfrac{\alpha}{n}x\big|\ge M\ge r+\diam(\Omega),
	\qquad\mbox{for any $x\in \partial \widetilde\Omega_{b,w}$.}
\end{equation}
Hence, Proposition~\ref{propertyf}\,(3) 
ensures that $f^\ast$ is of class $C^{1,1}$ in a neighborhood of $\frac{\alpha}{n}x$ whenever $x\in \partial \widetilde\Omega_{b,w}$. 
Next, we compute the principal curvatures of $\partial\widetilde\Omega_{b,w}$. To this aim we need to know that for any $x\in\partial\widetilde\Omega_{b,w}$ there holds
$$
    \nabla f^\ast \Big(\frac{\alpha}{n}\,x\Big)
    - w \not= 0.
$$
Indeed, due to \eqref{boundary-point} and \eqref{nfs}, we have
\begin{equation}\label{est-grad}
    \min_{x\in\partial\widetilde\Omega_{b,w}} \Big|\nabla f^\ast \Big(\frac{\alpha}{n}\,x\Big)
    - w\Big|>\frac{R}{\epsilon}.
\end{equation}
Consequently, for $x\in\partial\widetilde\Omega_{b,w}$ the outward  pointing normal $\widetilde\nu_x$ to $\partial\widetilde\Omega_{b,w}$ in $x$ is given by
$$
   \widetilde\nu_x =\frac{\nabla f^\ast \big(\frac{\alpha}{n}x\big)
    - w}{\big|\nabla f^\ast \big(\frac{\alpha}{n}x\big)
    - w\big|}, 
$$
and the tangent space of $\partial\widetilde\Omega_{b,w}$ at $x$ is given by $T_x\partial\widetilde\Omega_{b,w}=
(\widetilde\nu_x\R)^\perp$. 
The second fundamental form $\widetilde A_x\colon T_x\partial\widetilde\Omega_{b,w}\times T_x\partial\widetilde\Omega_{b,w}\to  \R$ of $\partial\widetilde\Omega_{b,w}$ at $x$ is given by 
\begin{equation}\label{second-fund}
    \widetilde A_x(\xi,\zeta) =-\frac{D^2f^\ast \big(\frac{\alpha}{n}x\big)(\xi,\zeta)}{\big|\nabla f^\ast \big(\frac{\alpha}{n}x\big)
    - w\big|},
    \qquad\forall\,\xi,\eta\in T_x\partial\widetilde\Omega_{b,w},
\end{equation}
provided $D^2f^\ast\big(\frac{\alpha}{n} x\big) $ exists. This, however, cannot be guaranteed in general. 
Therefore, we consider $\widetilde\Omega_{\tau,w}$ with $\tau\in[b,b+\delta)$ instead of $\widetilde\Omega_{b,w}$ and show by a co-area formula type argument the existence of second derivatives $\mathcal H^{n-1}$-a.e.~on $\partial\widetilde\Omega_{\tau,w}$ for a.e.~$\tau\in [b,b+\delta)$. 
Indeed, using the fact that $\partial\widetilde\Omega_{\tau,w}$ is the level set
$h_{w}^{-1} \{ \tau\}$ of the function $h_{w}(x):= \frac{n}{\alpha} f^\ast \big(\frac{\alpha}{n}x\big) -w\cdot x$, the co-area formula implies for any $\delta>0$ that
\begin{align*}
    \int_{b}^{b+\delta} \mathcal H^{n-1}\big( h_{w}^{-1} \{ \tau\}\cap \Sigma\big)\, \dtau
    &=
    \int_{ (\widetilde\Omega_{w,b+\delta}\setminus \widetilde\Omega_{b,w})\cap \Sigma}|\nabla h_{w}(x)|\, \dx\\
    &=
    \int_{ (\widetilde\Omega_{w,b+\delta}\setminus \widetilde\Omega_{b,w})\cap \Sigma}\big|\nabla f^\ast \big(\tfrac{\alpha}{n}x\big)-w\big|\, \dx,
\end{align*} 
where $\Sigma$ denotes the set of points in which $D^2f^\ast(\frac{\alpha}{n}x)$ does not exist. This,
however, is a set of $\mathcal L^n$-measure zero by Rademacher's theorem, and therefore the right-hand side above is equal to zero. This proves
that $\mathcal H^{n-1}( h_{w}^{-1} \{ \tau\}\cap \Sigma)=0$ for a.e.~$\tau\in [b,b+\delta)$.
Therefore we can choose a sequence $(\tau_i)_{i\in\N}\subset [b,b+\delta)$ such that $\tau_i\downarrow b$ and 
$D^2f^\ast(\frac{\alpha}{n}x)$ exists for $\mathcal H^{n-1}$-almost every $x\in \partial \widetilde \Omega_{\tau_i, w}$. Due to the continuity of $f^\ast$ and $\nabla f^\ast$ in a neighborhood of $\frac{\alpha}{n}(x_o-y(t))$ there exist $x_i\in \R^n$ and $\nu_i\in\R^n$ such that $x_i\in y(t)+\partial \widetilde\Omega_{\tau_i, w}$, $\nu_i$ is the outward pointing normal of $y(t)+\partial \widetilde\Omega_{\tau_i, w}$ at $x_i$, $x_i\to x_o$ and $\nu_i\to\nu_{x_o}$ as $i\to\infty$. 
Using Proposition \ref{propertyf}\,(3) and \eqref{est-grad} to estimate the right-hand side in \eqref{second-fund}, we obtain
for $\mathcal{H}^{n-1}$-a.e.~$x\in\partial\widetilde\Omega_{\tau_i,w}$ that
\begin{align}\label{curvature}
    \widetilde A_x(\xi,\xi) 
    &=-
    \frac{D^2f^\ast \big(\frac{\alpha}{n}x\big)(\xi,\xi)}{\big|\nabla f^\ast \big(\frac{\alpha}{n}x\big)- w\big|}
    \le
    \frac{\frac1{\epsilon}|\xi|^2}{\frac{R}{\epsilon}}
    =
    \tfrac{1}{R}|\xi|^2
    \quad\forall\, \xi\in\R^n. 
\end{align}
This ensures that for $\mathcal{H}^{n-1}$-a.e.~$x\in\partial\widetilde\Omega_{\tau_i,w}$
the principal curvatures of $\partial\widetilde\Omega_{\tau_i,w}$ are smaller than $\frac1R$, which is the principal curvature of $B_R(x_i-R\nu_{i})$. By the convexity of $\widetilde\Omega_{\tau_i,w}$ and the fact that $f^\ast$ is of class $C^{1,1}$ this implies $B_R(x_i-R\nu_{i})\subset y(t)+\widetilde\Omega_{\tau_i,w}$ for any $i\in\N$. Passing to the limit $i\to\infty$, we conclude $B_R(x_o-R\nu_{x_o})\subset y(t)+\widetilde\Omega_{b,w}=\widetilde\Omega_t$, which in view of Remark~\ref{rem:unif-conv} implies the set inclusion~\eqref{step3}. 
%
%
%
\medskip

\textit{Step 4: Proof of claim~(3).} 
In this step we will choose $\alpha_o$ in such a way that $v$ is a sub-solution to~\eqref{eq:strong-Cauchy-Dirichlet-prob} whenever $\alpha\ge\alpha_0$. 
We compute 
\begin{align*}
    \partial_t v(x,t)
    &=
    \partial_t 
    \Big[\frac{n}{\alpha}f^*\left(\frac{\alpha}{n}\big(x-y(t)\big)\right)-c(t)\Big] \\
    &=
    -\Big[\nabla f^*\left(\frac{\alpha}{n}\big(x-y(t)\big)\right)-\nabla f^*\left(\frac{\alpha}{n}\big(x_o-y(t)\big)\right)\Big]\cdot y'(t) +
    \partial_t g(x_o,t).
\end{align*}
From~\eqref{zlambdaesplicit} and the choices of $\rho_o$ and $M$ we know that 
\begin{equation*}
    |x_o-y(t)|
    =
    |z_\lambda(t)|
    \ge 
    \rho_o
    \ge 
    M
    \ge 
    r+\diam(\Omega).
\end{equation*} 
Since $\alpha\ge \alpha_o\ge n$ this ensures $\frac{\alpha}{n}(x-y(t))\in\R^n\setminus B_r$ for any $x\in\Omega$. Since $\Omega$ is convex, we thus have that $\frac{\alpha}{n}(\xi-y(t))\in\R^n\setminus B_r$ for any $\xi\in[x,x_o]$. 
Therefore, we may use Proposition~\ref{propertyf}\,(3) to conclude that
\begin{align*}
    |\partial_t v(x,t)|
    &\le 
    \frac{1}{\epsilon} \frac{\alpha}{n}|x-x_o| |y'(t)| +
    \|\partial_t g\|_{L^\infty(\Omega_T)} \\
    &= 
    \frac{1}{\epsilon} \frac{\alpha}{n}|x-x_o| |z_\lambda'(t)| +
    \|\partial_t g\|_{L^\infty(\Omega_T)} \\
    &\le 
    \frac{1}{\epsilon} \diam(\Omega) \big|D^2f\big(\widetilde w_{x_o}^-(t)+\lambda\nu_{x_o}\big)\big| |(\widetilde w_{x_o}^-)'(t)| +
    \|\partial_t g\|_{L^\infty(\Omega_T)} .
\end{align*}
Note that the choice of $\lambda$ ensures that $|\widetilde w_{x_o}^-(t)+\lambda\nu_{x_o}|\ge \lambda - Q_1\ge 1$ and hence $\widetilde w_{x_o}^-(t)+\lambda\nu_{x_o}\in \R^n\setminus B_1$. 
Since $D^2f$ is bounded outside the unit ball by assumption \ref{A2} and $|(\widetilde w_{x_o}^-)'(t)|\le \mathsf Q$ by \ref{A3} we finally conclude
\begin{align*}
    |\partial_t v(x,t)|
    &\le 
    \frac{1}{\epsilon}\diam(\Omega) \|D^2f\|_{L^\infty(\R^n\setminus B_1)}\mathsf Q +
    \|\partial_t g\|_{L^\infty(\Omega_T)}.
\end{align*}
We choose
\begin{equation}\label{choice-alpha-o}
	\alpha_o
	:=
	\max\bigg\{n \,, \, 
	\frac{1}{\epsilon} \diam(\Omega) \|D^2f\|_{L^\infty(\R^n\setminus B_1)}  \mathsf Q+
    \|\partial_t g\|_{L^\infty(\Omega_T)} \bigg\}.
\end{equation}
Note that $\alpha_o$ depends on $n,\epsilon, \diam(\Omega), \|D^2f\|_{L^\infty(\R^n\setminus B_1)},  \|\partial_t g\|_{L^\infty(\Omega_T)}$ and $\mathsf Q$. 
In view of the preceding computation and \eqref{cellina} we have for any $\alpha\ge\alpha_o$ that 
\begin{equation}
    \partial_t v(x,t)-\Div\nabla f\big(\nabla v(x,t)\big)\le |\partial_t v(x,t)|-\alpha\le 0 ,
\end{equation}
i.e.~$v$ is sub-solution of the parabolic equation~\eqref{eq:strong-Cauchy-Dirichlet-prob} and hence claim (3) is proved. 
\medskip 

\textit{Step 5: Proof of claim~(4).}  With~\eqref{def:y(t)} and~\eqref{zlambdaimplicit} we compute
\begin{align*}
   \nabla v(x,t)
    &=
    \nabla f^*\left(\frac{\alpha}{n}\big(x-y(t)\big)\right)\\
    &=
    \nabla f^*\left(\frac{\alpha}{n}\big(x-y(t)\big)\right)
    -
    \nabla f^*\left(\frac{\alpha}{n}\big(x_o-y(t)\big)\right)
    +
    \nabla f^*\Big(\frac{\alpha}{n}z_\lambda (t)\Big) \\
    &=
    \nabla f^*\left(\frac{\alpha}{n}\big(x-y(t)\big)\right)
    -
    \nabla f^*\left(\frac{\alpha}{n}\big(x_o-y(t)\big)\right) +
    \lambda\nu_{x_o} + 
    \widetilde w_{x_o}^-(t).
\end{align*}
The difference of the first two terms is bounded exactly as in Step~4, so that 
\begin{align*}
    |\nabla v(x,t)|
    &\le
    \frac{\alpha}{\epsilon n}\diam (\Omega) +
    \lambda +\| \widetilde w_{x_o}^{-}\|_{L^\infty ([0,T],\R^n)} \\
     &\le
    \frac{\alpha}{\epsilon n}\diam (\Omega) +
    \lambda + Q_1.
\end{align*}
 This implies
\begin{align*}
    |\nabla v(x,t)|
    &\le c\big(n, \epsilon, R, \diam(\Omega), f, \nabla f, Q_1, \alpha\big).
\end{align*}
This proves claim (4) and finishes the proof of Proposition~\ref{prop:barrier}.
\end{proof}

\subsection{Parabolic sub- super-minimizers and the comparison principle}\label{sec:uniqueness}

We recall that the variational solution constructed in Proposition~\ref{prop:existence-reg} admits a time derivative $\partial_tu\in L^2(\Omega_T)$. Therefore, we may perform an integration by parts in the first term on the right-hand side of the variational inequality \eqref{eq:variational-inequality-rest}. In this way, the variational inequality can be re-written as
\begin{align*}
	\iint_{\Omega_T}f(\nabla u)\,\dx\dt
	& \le
	\iint_{\Omega_T}\big[\partial_{t}u(v-u)
	+
	f(\nabla v)\big]\,\dx\dt
\end{align*}
for any $v\in \Lip_g(\Omega_T, \LL)$ with $\partial_{t}v\in L^{2}(\Omega_{T})$. 
Note that the assumption $\partial_{t}v\in L^{2}(\Omega_{T})$ can be eliminated by an approximation argument. This motivates  the following definition:

\begin{definition}[Parabolic sub-/super-minimizer]\label{def:min-rest}\upshape
Let $\LL\in(0,\infty]$. A map
$u\in \Lip(\Omega_T, \LL)$ with $\partial_t u\in L^2(\Omega_T)$ is called \emph{parabolic sub-minimizer
(of the gradient constrained obstacle problem in the case $L<\infty$)} in $\Lip(\Omega_T, \LL)$ if and only if 
\begin{align}\label{eq:minimizer-rest}
	\iint_{\Omega_T}f(\nabla u)\,\dx\dt
	& \le
	\iint_{\Omega_T}\big[\partial_{t}u(v-u) +
	f(\nabla v)\big]\,\dx\dt
\end{align}
holds true for any  $v\in \Lip_u(\Omega_T, \LL)$ with $v\le u$ in $\Omega_T$.
Moreover, a map
$u\in \Lip(\Omega_T, \LL)$ is called \emph{parabolic super-minimizer} in $\Lip(\Omega_T, \LL)$ if and only if \eqref{eq:minimizer-rest} holds true for any  $v\in \Lip_u(\Omega_T, \LL)$ with $v\ge u$ in $\Omega_T$. Finally, $u\in \Lip(\Omega_T, \LL)$ is called \emph{parabolic minimizer} if and only if \eqref{eq:minimizer-rest} holds true for any  $v\in \Lip_u(\Omega_T, \LL)$.
\end{definition}

The concept of parabolic minimizers for vector-valued integrands with quadratic growth originates from the work of Wieser \cite{Wieser.1987}. Subsequently, it will be essential to establish that a localization principle with respect to the spatial variables holds for parabolic minimizers

\begin{remark}[Localization in space]\upshape\label{rem:locx-1}
Let $\LL\in(0,\infty)$ and suppose that $\Omega'\subset\Omega$ is an open, convex set, $f\colon\R^n\to \R$ a convex integrand and $u\in \Lip(\Omega_T, \LL)$ with $\partial_t u\in L^2(\Omega_T)$ a parabolic minimizer in the sense of Definition \ref{def:min-rest} in $\Omega_T$. Then, $u$ is also a parabolic minimizer in the subcylinder $\Omega'_T:=\Omega'\times[0,T)$. 
The proof of this elementary fact is analogous to the setting of time independent boundary data and can be found in \cite[Remark~4.2]{BDMS_bd-slope}. 
\end{remark}

In the following lemma, we establish the comparison principle for parabolic sub- and super-minimizers.

\begin{lem}[Comparison principle]\label{lem:comp-rest}
Let $\LL\in(0,\infty]$ and suppose that $\Omega\subset\R^n$ is a bounded open set, $f\colon\R^n\to \R$ is a convex integrand and $u, \tilde u\in \Lip(\Omega_T, \LL)$ with $\partial_t u,\partial_t\tilde u\in L^2(\Omega_T)$. Suppose that $u$ is a sub-minimizer and $\tilde u$ is a super-minimizer in $\Omega_T$ in the sense of Definition \ref{def:min-rest}  and that $u\le \tilde u$ on $\partial_{\mathcal P}\Omega_T$. Then, we have
$$
	u\le \tilde u\quad\mbox{a.e. in $\Omega_T$.}
$$
\end{lem}

\begin{proof}[Proof]
Let $\tau\in(0,T]$. 
We define 
\begin{equation*}
    v:=
    \left\{\begin{array}{cl}
        \min\{u, \tilde u\}, &
        \mbox{in $\Omega_\tau$,} \\[5pt]
        u, &
        \mbox{in $\Omega\times[\tau,T)$,}
    \end{array}\right.
\end{equation*}
and 
\begin{equation*}
    w:=
    \left\{\begin{array}{cl}
        \max\{u, \tilde u\}, &
        \mbox{in $\Omega_\tau$,} \\[5pt]
        \tilde u, &
        \mbox{in $\Omega\times[\tau,T)$.}
    \end{array}\right.
\end{equation*}
We note that $v\in \Lip_u(\Omega_T, \LL)$ with $v\le u$ and $v(0)=u_o$ and $w\in \Lip_{\tilde u}(\Omega_T, \LL)$ with $w\ge\tilde u$ and $w(0)=\tilde u_o$.
This ensures that $v$ is an admissible comparison function in the variational inequality  \eqref{eq:minimizer-rest} for $u$ and $w$ for  of $\tilde u$. Adding the two resulting inequalities and taking into account that the parts of the integrals on $\Omega\times(\tau,T)$ cancel themselves out, we obtain
\begin{align}\label{min-comp}
	&\iint_{\Omega_\tau}\big[f(\nabla u) + f(\nabla \tilde u)\big] \,\dx\dt \nonumber\\
	&\qquad\le
	\iint_{\Omega_\tau}\big[f(\nabla v) + f(\nabla w) +
	\partial_t u (v-u) + \partial_t \tilde u(w-\tilde u)\big]\,\dx \dt.
\end{align}
We now consider the terms on the right-hand side of \eqref{min-comp}. From the definition of $v$ and $w$ we infer 
\begin{align*}
	\iint_{\Omega_\tau}\big[f(\nabla v) + f(\nabla w)\big] \,\dx\dt
	= 
	\iint_{\Omega_\tau}\big[f(\nabla u) + f(\nabla \tilde u)\big] \,\dx\dt.
\end{align*}
Moreover, we observe that $v-u=-(u-\tilde u)_+$ and $w-\tilde u=(u-\tilde u)_+$ in $\Omega_\tau$, so that 
\begin{align*}
	\partial_t u(v-u) + \partial_t \tilde u(w-\tilde u)
	&=
	-\partial_t (u-\tilde u)(u-\tilde u)_+ 
    =
	-\tfrac12\partial_t (u-\tilde u)_+^2 .
\end{align*}
This implies 
\begin{align*}
	\iint_{\Omega_\tau} 
	\big[\partial_t u(v-u) + \partial_t \tilde u(w-\tilde u)\big] \dx\dt 
	&=
	-\tfrac12\iint_{\Omega_\tau}
	\partial_t (u-\tilde u)_+^2 \dx\dt \\
	&=
	-\tfrac12 \int_{\Omega\times\{\tau\}} 
	(u-\tilde u)_+^2 \dx .
\end{align*}
Here we used the assumption that $(u-\tilde u)_+(0)=0$.
Joining the preceding identities with \eqref{min-comp}, we conclude 
\begin{align*}
	\int_{\Omega\times\{\tau\}}(u - \tilde u)_+^2\dx 
	\le
	0 .
\end{align*}
Since $\tau\in(0,T]$ was arbitrary, this proves the claim $u\le \tilde u$ a.e.~in $\Omega_T$.
\end{proof}

As a consequence of the preceding comparison principle, we obtain the following result.

\begin{lem}[Maximum principle]\label{lem:maxpr}
Let $\LL\in(0,\infty]$ and suppose that $\Omega\subset\R^n$ is open and bounded, $f\colon\R^n\to \R$ convex and let $u, \tilde u\in \Lip(\Omega_T, \LL)$ with $\partial_t u,\partial_t\tilde u\in L^2(\Omega_T)$. Suppose that $u$ is a sub-minimizer and $\tilde u$ is a super-minimizer in the sense of Definition \ref{def:min-rest} in $\Omega_T$. Then, we have
$$
	\sup_{\Omega_T}(u-\tilde u)= \sup_{\partial_{\mathcal{P}}\Omega_{T}}(u-\tilde u).
$$
\end{lem}

\begin{proof}[Proof]
For $(x,t)\in\partial_{\mathcal P}\Omega_T$ there holds 
$$
	u(x,t)
	=
	\tilde u(x,t) + u(x,t) - \tilde u(x,t)
	\le 
	\tilde u(x,t) + \sup_{\partial_{\mathcal P}\Omega_T}(u - \tilde u).
$$
Since $\tilde u$ is a parabolic super-minimizer in $\Omega_T$ it immediately follows  that also $\tilde u+\sup_{\partial_{\mathcal P}\Omega_T}(u - \tilde u)$ is parabolic super-minimizer in $\Omega_T$.
In view of Lemma \ref{lem:comp-rest} we therefore have 
$$
	u
	\le 
	\tilde u+\sup_{\partial_{\mathcal P}\Omega_T}(u - \tilde u)
	\quad\mbox{in $\Omega_T$.}
$$
This inequality can be re-written in the form 
$$
	\sup_{\Omega_T}(u - \tilde u)
	\le 
	\sup_{\partial_{\mathcal P}\Omega_T}(u - \tilde u).
$$
Since the reversed inequality holds trivially, this proves the claim.
\end{proof}

\subsection{A quantitative bound on the Lipschitz continuity}\label{sec:lip}

In this subsection, we establish a quantitative Lipschitz bound for parabolic minimizers of the gradient constraint problem under the 
$t-{\rm BSC}_Q$.

\begin{prop}\label{prop:Lip}
Let the  assumptions of Theorem~\ref{thm:existence-var-sol} be in force and $L\in(0,\infty]$. Then, every parabolic minimizer $u\in \Lip(\Omega_T, \LL)$ with $\partial_t u\in L^2(\Omega_T)$ in the sense of Definition~\ref{def:min-rest} and Cauchy-Dirichlet boundary datum $g$ satisfies the  gradient bound
$$
	\|\nabla u\|_{L^\infty (\Omega_T,\R^n)}
	\le 
	C,
$$
where $C$ depends on 
$n, \epsilon, R, \diam(\Omega), f, \nabla f, \|D^2f\|_{L^\infty(\R^n\setminus B_1)}, Q, [g]_{0,1;\Omega_T}, \mathsf Q$
\end{prop}

\begin{proof}[Proof]
Let $x_1\not= x_2$ two arbitrary points in $\Omega$ and $t\in (0,T)$. Define $y:=x_2-x_1$ and
\begin{align*}
	u_y(x,t)
	:=
	u(x+y,t),
	\qquad\mbox{for $(x,t)\in \widetilde\Omega_T$},
\end{align*}
where $\widetilde\Omega_T:=\{(x-y,t)\in\R^{n+1}: (x,t)\in\Omega_T\}$. 
Then, $u_y$ is a parabolic minimizer in $\widetilde\Omega_T$ in the class $\Lip(\widetilde\Omega_T, \LL)$ in the sense of Definition~\ref{def:min-rest}.
We denote the intersection of both cylinders by $(\Omega\cap\widetilde\Omega)_T:=(\Omega\cap\widetilde\Omega)\times(0,T)$. 
Lemma~\ref{rem:locx-1} ensures that both, $u$ and $u_y$ are parabolic minimizers in $(\Omega\cap\widetilde\Omega)_T$ in $\Lip((\Omega\cap\widetilde\Omega)_T, \LL)$. Therefore, from the maximum principle in Lemma \ref{lem:maxpr} we conclude that there exists a boundary point $(x_o,t_o)\in \partial_{\mathcal P}\big((\Omega\cap\widetilde\Omega)_T\big)$ such that 
\begin{align*}
	|u(x_1,t)-u_y(x_1,t)|
	\le
	|u(x_o,t_o)-u_y(x_o,t_o)|.
\end{align*}
In view of the definition of $u_y$, this inequality yields
\begin{align}\label{bound-diff-u}
	|u(x_1,t)-u(x_2,t)|
	\le
	|u(x_o,t_o)-u(x_o+y,t_o)|.
\end{align}
Since $(x_o,t_o)\in \partial_{\mathcal P}\big((\Omega\cap\widetilde\Omega)_T\big)$, we either have $t_o=0$ or $x_o\in \partial (\Omega\cap\widetilde\Omega)$.
In the first case, we recall that $u(\cdot,0)=g(\cdot,0)$.  Using the Lipschitz condition of $g$, we then obtain
\begin{align*}
	\big|u(x_1,t)-u(x_2,t)\big|
	\le
	|g(x_o,0)-g(x_o+y,0)| 
	\le
	\|\nabla g_o\|_{L^\infty (\Omega,\R^n)}|y|.
\end{align*}
In the other  case we know that one of the points $x_o$ or $x_o+y$ belongs to $\partial\Omega$. Without loss of generality we may assume $x_o\in \partial\Omega$. Since $u=g$ on the lateral boundary of $\Omega_T$, inequality \eqref{bound-diff-u} turns into \begin{align*}
	|u(x_1,t)-u(x_2,t)|
	\le
	|g(x_o,t_o)-u(x_o+y,t_o)|.
\end{align*}
By $v_{x_o}^\pm$ we denote the barrier functions constructed in Proposition~\ref{prop:barrier} applied with $\alpha=\pm\alpha_o$. The barrier functions satisfy $v_{x_o}^\pm(x_o,t)=g(x_o,t)$ for any $t\in[0,T)$ and $v_{x_o}^-\le g\le v_{x_o}^+$ in $\Omega_T$. Moreover, $v_{x_o}^-$ is a sub-solution and since $f$ is convex, it is also a parabolic sub-minimizer in the sense of Definition~\ref{def:min-rest}. Similarly, $v_{x_o}^+$ is a super-minimizer.
Furthermore, $v_{x_o}^\pm$ are Lipschitz continuous with respect to the spatial variable
with Lipschitz constant $\widetilde Q$ depending on 
$n, \epsilon, R, \diam(\Omega), f, \nabla f$, $\|D^2f\|_{L^\infty(\R^n\setminus B_1)}, Q, [g]_{0,1;\Omega_T}, \mathsf Q$.
Then, the maximum principle from Lemma~\ref{lem:maxpr} implies
\begin{equation*}
	v_{x_o}^-\le u\le v_{x_o}^+\qquad \mbox{in $\Omega_T$,}
\end{equation*}
so that
\begin{align*}
    -\widetilde Q|y|
    &\le v_{x_o}^+(x_o,t)- v_{x_o}^+(x_o+y,t) \\
    &\le
     g(x_o,t)- u(x_o+y,t)\\
    &\le 
    v_{x_o}^-(x_o,t)- v_{x_o}^-(x_o+y,t)
    \le 
    \widetilde Q|y|.
\end{align*}
Therefore,  we have
\begin{align*}
	|u(x_1,t)-u(x_2,t)|
	\le
	\widetilde Q|y|.
\end{align*}
Joining both cases and recalling that $y=x_2-x_1$, we obtain
$$
	|u(x_1,t)-u(x_2,t)|
	\le
	\max\big\{ \widetilde Q, \| \nabla g_o\|_{L^\infty (\Omega,\R^n)}\big\} |x_1-x_2|,
$$
which implies the claimed gradient bound.
\end{proof}

\subsection{Proof of Theorem~\ref{thm:existence-var-sol}}
In this section, we indicate how the gradient hypothesis 
$\|\nabla u\|_{L^\infty (\Omega,\R^n)}\le \LL$ for a variational solution $u$
of the gradient-constrained obstacle problem in 
$\Lip_{g}(\Omega_{T}, \LL)$ can be removed, thereby obtaining the existence result from Theorem~\ref{thm:existence-var-sol}.


\begin{proof}[Proof of Theorem~\ref{thm:existence-var-sol}]
Let $\LL>C$, where $C$ denotes the constant from Proposition~\ref{prop:Lip} depending only on $n, \epsilon, R, \diam(\Omega), f, \nabla f, \|D^2f\|_{L^\infty(\R^n\setminus B_1)}, Q, [g]_{0,1;\Omega_T}, \mathsf Q$.
Due to Proposition~\ref{prop:ex-constr} there exists a unique variational solution $u\in \Lip_g(\Omega_T,\LL)$ of the gradient constrained problem.
The solution admits a weak time derivative $\partial_tu\in L^2(\Omega_T)$ and satisfies $u(0)=g_o$ in the $L^2(\Omega)$-sense. As pointed out at the beginning of \S\,\ref{sec:uniqueness},  $u$ is also a parabolic minimizer in the sense of Definition~\ref{def:min-rest}. Proposition~\ref{prop:Lip} ensures  that the strict gradient bound 
$$
	\|\nabla u\|_{L^\infty (\Omega_T,\R^n)}
	\le 
	C
	<
	L
$$
holds. Therefore, it
remains to prove that the variational inequality \eqref{eq:variational-inequality-rest} satisfied by $u$ actually holds for any comparison map $w\in \Lip_{g}(\Omega_{T})$ with $\partial_t w\in L^2(\Omega_T)$.  
To this aim we consider $w\in \Lip_{g}(\Omega_{T})$ and define
$$
	v:= u+s(w-u)\qquad \mbox{for $0<s\ll 1$.}
$$
Observe that $v$ is an admissible comparison function in \eqref{eq:minimizer-rest}, since $v$ coincides with $u$ on the lateral boundary $\partial\Omega\times(0,T)$ and $\|\nabla v\|_{L^\infty (\Omega_T,\R^n)}<\LL$ for $s>0$
small enough. From \eqref{eq:minimizer-rest} and the convexity of $f$ we infer 
\begin{align*}
	\iint_{\Omega_T}f(\nabla u)\,\dx\dt
	& \le
	\iint_{\Omega_T}\big[s\partial_{t}u(w-u)
	+
	f\big( (1-s)\nabla u+s\nabla w\big)\big]\,\dx\dt \\
	& \le
	\iint_{\Omega_T}\big[s\partial_{t}u(w-u)
	+
	(1-s)
	f(\nabla u)+sf(\nabla w)\big]\,\dx\dt .
\end{align*}
We re-absorb the second term of the right-hand side into the left
and divide the result by $s>0$, so that 
\begin{align*}
	\iint_{\Omega_T}f(\nabla u)\,\dx\dt
	& \le
	\iint_{\Omega_T}\big[\partial_{t}u(w-u)
	+
	f(Dw)\big]\,\dx\dt \\
 	&=
	\iint_{\Omega_T}\big[\partial_{t}w(w-u)-\tfrac12\partial_t|w-u|^2
	+
	f(\nabla w)\big]\,\dx\dt \\
	& =
	\iint_{\Omega_T}\big[\partial_{t}w(w-u)
	+
	f(\nabla w)\big]\,\dx\dt \\
	& \quad+
	\tfrac{1}{2}\|w(0)-g_{o}\|_{L^{2}(\Omega)}^{2}
	-
	\tfrac{1}{2}\|(w-u)(T)\|_{L^{2}(\Omega)}^{2}.
\end{align*}
This shows that the variational inequality \eqref{eq:variational-inequality-rest} holds for every comparison function $w\in \Lip_{g}(\Omega_{T})$ that satisfies $\partial_tw\in L^2(\Omega_T)$. Therefore, $u$ is a variational solution of the unconstrained problem in the sense of Definition~\ref{def:var-sol}. As before, the uniqueness can be deduced as in~\cite[Lemma~3.3]{BDS} and the proof of Theorem~\ref{thm:existence-var-sol} is complete.
\end{proof}

\section{Regularity of solutions and proof of Theorem~\ref{thm:reg-gen}}\label{sec:reg}


Our aim in this section is to prove Theorem~\ref{thm:reg-gen}. By $u$ we denote the unique variational solution from Theorem~\ref{thm:existence-var-sol} satisfying
$$
	\|\nabla u\|_{L^\infty (\Omega_T,\R^n)}
	\le 
	C
    =:
    M.
$$
Since $f$ is assumed to be of class $C^1$ and $\nabla u$ is bounded, the associated Euler-Lagrange equation is well defined. Hence, $u$ is a  weak solution of the parabolic Cauchy-Dirichlet problem~\eqref{eq:strong-Cauchy-Dirichlet-prob}. 
From the Poincar\'e inequality for solutions to parabolic equations (cf.~\cite[Lemma 3.1]{Boegelein-Duzaar-Mingione}), we obtain
\begin{align*}
	\biint_{Q_\varrho(z_o)}|u-(u)_{z_o;\varrho}|^2\dx\dt
	&\le
	C(n)\,\varrho^2\bigg[
	\biint_{Q_\varrho(z_o)}|\nabla u|^2\dx\dt
	+\sup_{B_M (0)}|Df|^2\bigg]
	\le
	C\varrho^{2},
\end{align*}
for any parabolic cylinder $Q_\varrho(z_o):=B_\varrho(x_o)\times(t_o-\varrho^2,t_o+\rho^2)\subset \Omega_T$
with $z_o=(x_o,t_o)\in\Omega_T$. A similar Poincar\'e inequality holds for parabolic cylinders with center $(x_o,t_o)\in\partial_{\mathcal P}\Omega_T$. This can be deduced as in \cite{BDM1}, since $\Omega$ is convex and therefore also a Lipschitz domain. 
The parabolic version of Campanato's characterization of H\"older continuity (with respect to the parabolic metric) by Da Prato \cite{Prato} implies the Lipschitz continuity of $u$ with respect to the parabolic metric, i.e.~$u\in C^{0;1,\frac12}(\overline\Omega_T)$. This finishes the proof of Theorem~\ref{thm:reg-gen}.

\end{document}